\documentclass[a4paper,pdftex,12pt]{article}
\usepackage{times}
\usepackage{amsmath}
\usepackage{graphicx}
\usepackage{caption}
\usepackage{subcaption}
\usepackage{tikz}
\usepackage{hyperref}
\hypersetup{linktocpage,colorlinks=true,pdfborder={0 0 0},linkcolor=black,urlcolor=black}
\usepackage{geometry}
\title{New minimal $(4;n)$-regular matchstick graphs}

\author{Mike Winkler$^\star$\quad Peter Dinkelacker$^{\star\star}$\quad Stefan Vogel$^{\star\star\star}$}

\date{\small May 1, 2018\\[4mm]$^\star$Fakult\"at f\"ur Mathematik, Ruhr-Universit\"at Bochum, Germany, mike.winkler@ruhr-uni-bochum.de\\[2mm]
            $^{\star\star}$Togostr. 79, 13351 Berlin, Germany, peter@grity.de\\[2mm]
            $^{\star\star\star}$Raun, Dorfstr. 7, 08648 Bad Brambach, Germany, backebackekuchen16@gmail.com}

\begin{document}
  
  \maketitle
  
  {\small\noindent\textbf{Abstract:} A matchstick graph is a graph drawn with straight edges in the plane such that the edges have unit length, and non-adjacent edges do not intersect. We call a matchstick graph $(m;n)$-regular if every vertex has only degree $m$ or $n$. In this article the authors present the latest known $(4;n)$-regular matchstick graphs for $4\leq n\leq11$ with a minimum number of vertices.}
  
  \section{\large{Introduction}}
  \noindent
  A matchstick graph is a planar unit-distance graph. That is a graph drawn with straight edges in the plane such that the edges have unit length, and non-adjacent edges do not intersect. We call a matchstick graph $(m;n)$-regular if every vertex has only degree $m$ or $n$.
  
  For $m\leq n$ minimal $(4;n)$-regular matchstick graphs with a minimum number of vertices only exist for $4\leq n\leq11$. The smallest known $(4;n)$-regular matchstick graph for $n=4$, also named 4-regular, is the so-called \textit{Harborth graph} (Fig. 2) consisting of 52 vertices and 104 edges. Since 1986 the second smallest known 4-regular matchstick graph consists of 60 vertices and 120 edges (Fig. 5). On July 3, 2016 the authors presented a new second smallest known example with 54 vertices and 108 edges (Fig. 3), which is based on the Harborth graph [8]. On April 15, 2016 Mike Winkler presented an example with 57 vertices and 114 edges with a whole new geometry (Fig. 4) [9][10]. For each $n>11$ only infinite graphs with an infinite number of vertices exists. Figures 22 -- 25 show examples of infinite graphs for $n=12$ and $n=13$.
  
  It is an open problem how many different $(4;n)$-regular matchstick graphs with a minimum number of vertices for $4\leq n\leq11$ exist and which is the least minimal number. "Our knowledge on matchstick graphs is still very limited. It seems  to be hard to obtain rigid mathematical results about them. Matchstick problems constructing the minimal example can be quite challenging. But the really hard task is to rigidly prove that no smaller example can exist." [3]
  
  In this article the authors present the $(4;n)$-regular matchstick graphs for $5\leq n\leq11$ with the smallest currently known number of vertices as well as further new minimal examples for $n=4$, $n=5$ and $n=6$ with 108, 114, 118, 120, 121 and 126 edges. The earlier versions of these graphs were presented on the \textit{Math Magic}\footnote{http://www2.stetson.edu/~efriedma/mathmagic/1205.html} website of Erich Friedman [1], which can be seen on an older version of this site from 2015 [4].
  
  The authors discovered the new minimal $(4;n)$-regular matchstick graphs in the days from March 14, 2016 -- April 15, 2018 and presented them for the first time in a German mathematics internet forum [7].
  
  The geometry, rigidity or flexibility of the graphs in this article has been verified by Stefan Vogel with a computer algebra system named \textsc{Matchstick Graphs Calculator} (MGC) [6]. This remarkable software created by Vogel runs directly in web browsers. A special version of the MGC contains all graphs from this article and is available under this {weblink}\footnote{http://mikewinkler.co.nf/matchstick\_graphs\_calculator.htm}. The method Vogel used for the calculations he describes in a separate German article [5].
  
  \textit{Remark}: The MGC contains a constructive proof for each graph shown in this article. We are using this online reference, because these proofs are too extensive to reproduce here.
  
  \textit{Note}: In the PDF version of this article the vector graphics can be viewed with the highest zoom factor to see the smallest details. For example the very small rhombus in Figure 14.
  
  \section{\large{Rigid subgraphs}}
  
  The geometry of the new graphs, except for $n=11$, is not complicated. Most of these graphs have a high degree of symmetry. There are two types of most used rigid subgraphs, which we call the \textit{kite} (Fig. 1a) and the \textit{triplet kite} (Fig. 1d). The kite is a $(2;4)$-regular matchstick graph consisting of 12 vertices and 21 edges and has a vertical symmetry. Two kites can be connected to each other in two useful ways. We call these subgraphs the \textit{double kite} (Fig. 1b) and the \textit{reverse double kite} (Fig. 1c), both consisting of 22 vertices and 42 edges. The reverse double kite offers the possibility to connect two of the inner vertices with an additional unit length edge. This property has been used for $n=5$ in Figure 19 and slightly modified  in Figure 15. What makes the subgraphs (b) and (c) so useful is the fact that they have only two vertices of degree 2. Two of these subgraphs can be used to connect two vertices of degree 2 at different distances by using them like clasps. This property has been used for $n=9$, $n=10$ and $n=11$ (Fig. 11 -- 13).
  
  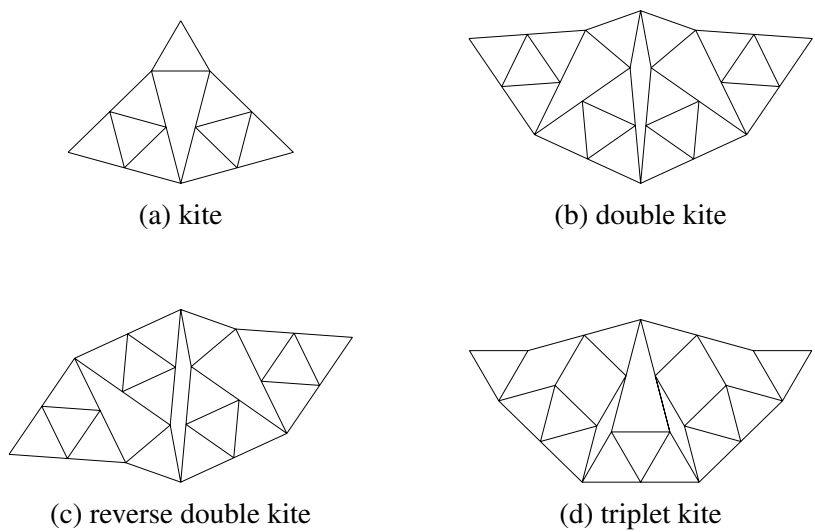
\begin{figure}[!ht]
    \centering
    \begin{minipage}[t]{0.45\linewidth}
      \centering
      \begin{tikzpicture}
      [y=0.5pt, x=0.5pt, yscale=-1.0, xscale=1.0]
      \draw[line width=0.01pt]
      (86.5346,1050.0456) -- (128.7085,1038.3319)
      (139.6510,995.9514) -- (128.7085,1038.3319)
      (128.7085,1038.3319) -- (97.4772,1007.6651)
      (97.4772,1007.6651) -- (86.5346,1050.0456)
      (97.4772,1007.6651) -- (139.6510,995.9514)
      (139.6510,995.9514) -- (108.4198,965.2847)
      (108.4198,965.2847) -- (97.4772,1007.6651)
      (108.4198,965.2847) -- (64.6494,965.2847)
      (86.5346,927.3784) -- (64.6494,965.2847)
      (64.6494,965.2847) -- (86.5346,1050.0456)
      (86.5346,1050.0456) -- (2.1869,1026.6183)
      (2.1869,1026.6183) -- (64.6494,965.2847)
      (75.5920,1007.6651) -- (33.4182,995.9515)
      (33.4182,995.9515) -- (44.3608,1038.3320)
      (44.3608,1038.3320) -- (75.5920,1007.6651)
      (128.7085,1038.3319) -- (170.8824,1026.6183)
      (170.8824,1026.6183) -- (139.6510,995.9514)
      (108.4198,965.2847) -- (86.5346,927.3784);
      \end{tikzpicture}
      \subcaption{kite}
    \end{minipage}
    \quad
    \begin{minipage}[t]{0.45\linewidth}
      \centering
      \begin{tikzpicture}
      [y=0.5pt, x=0.5pt, yscale=-1.0, xscale=1.0]
      \draw[line width=0.01pt]
      (54.8439,1011.4230) -- (30.2950,975.1850)
      (49.4036,935.8060) -- (30.2950,975.1850)
      (30.2950,975.1850) -- (73.9525,972.0441)
      (73.9525,972.0441) -- (54.8439,1011.4230)
      (73.9525,972.0441) -- (49.4036,935.8060)
      (49.4036,935.8060) -- (93.0611,932.6651)
      (93.0611,932.6651) -- (73.9525,972.0441)
      (93.0611,932.6651) -- (126.4125,961.0117)
      (134.2857,917.9552) -- (126.4125,961.0117)
      (126.4125,961.0117) -- (54.8439,1011.4230)
      (54.8439,1011.4230) -- (134.2857,1048.1977)
      (134.2857,1048.1977) -- (126.4125,961.0117)
      (90.6282,986.2174) -- (130.3491,1004.6047)
      (130.3491,1004.6047) -- (94.5648,1029.8103)
      (94.5648,1029.8103) -- (90.6282,986.2174)
      (30.2950,975.1850) -- (5.7461,938.9469)
      (5.7461,938.9469) -- (49.4036,935.8060)
      (93.0611,932.6651) -- (134.2857,917.9552)
      (213.7275,1011.4230) -- (238.2764,975.1850)
      (219.1678,935.8060) -- (238.2764,975.1850)
      (238.2764,975.1850) -- (194.6189,972.0441)
      (194.6189,972.0441) -- (213.7275,1011.4230)
      (194.6189,972.0441) -- (219.1677,935.8060)
      (219.1678,935.8060) -- (175.5102,932.6651)
      (175.5102,932.6651) -- (194.6189,972.0441)
      (175.5102,932.6651) -- (142.1589,961.0117)
      (134.2857,917.9552) -- (142.1589,961.0117)
      (142.1589,961.0117) -- (213.7275,1011.4230)
      (213.7275,1011.4230) -- (134.2857,1048.1977)
      (134.2857,1048.1977) -- (142.1589,961.0117)
      (177.9432,986.2174) -- (138.2223,1004.6047)
      (138.2223,1004.6047) -- (174.0066,1029.8103)
      (174.0066,1029.8103) -- (177.9432,986.2174)
      (238.2764,975.1850) -- (262.8253,938.9469)
      (262.8253,938.9469) -- (219.1677,935.8060)
      (175.5102,932.6651) -- (134.2857,917.9552);
      \end{tikzpicture}
      \subcaption{double kite}
    \end{minipage}
    
    \begin{center}\end{center}
    
    \begin{minipage}[t]{0.45\linewidth}
      \centering
      \begin{tikzpicture}
      [y=0.5pt, x=0.5pt, yscale=-1.0, xscale=1.0]
      \draw[line width=0.01pt]
      (53.0582,956.1585) -- (92.7791,937.7712)
      (128.5634,962.9769) -- (92.7791,937.7712)
      (92.7791,937.7712) -- (88.8425,981.3641)
      (88.8425,981.3641) -- (53.0582,956.1585)
      (88.8425,981.3641) -- (128.5634,962.9769)
      (128.5634,962.9769) -- (124.6268,1006.5698)
      (124.6268,1006.5698) -- (88.8425,981.3641)
      (124.6268,1006.5698) -- (91.2754,1034.9164)
      (132.5000,1049.6262) -- (91.2754,1034.9164)
      (91.2754,1034.9164) -- (53.0582,956.1585)
      (53.0582,956.1585) -- (3.9604,1028.6346)
      (3.9604,1028.6346) -- (91.2754,1034.9164)
      (72.1668,995.5375) -- (47.6179,1031.7755)
      (47.6179,1031.7755) -- (28.5093,992.3965)
      (28.5093,992.3965) -- (72.1668,995.5375)
      (92.7791,937.7712) -- (132.5000,919.3839)
      (176.1575,987.6460) -- (211.9418,1012.8517)
      (176.1575,987.6460) -- (136.4366,1006.0333)
      (136.4366,1006.0333) -- (140.3732,962.4403)
      (140.3732,962.4403) -- (176.1575,987.6460)
      (140.3732,962.4403) -- (173.7245,934.0937)
      (140.3732,962.4403) -- (132.5000,919.3839)
      (132.5000,919.3839) -- (173.7245,934.0937)
      (173.7245,934.0937) -- (211.9418,1012.8517)
      (211.9418,1012.8517) -- (261.0396,940.3756)
      (261.0396,940.3756) -- (173.7246,934.0937)
      (192.8332,973.4727) -- (217.3821,937.2346)
      (217.3821,937.2346) -- (236.4907,976.6136)
      (236.4907,976.6136) -- (192.8332,973.4727)
      (211.9418,1012.8517) -- (172.2209,1031.2390)
      (172.2209,1031.2390) -- (176.1575,987.6460)
      (172.2209,1031.2390) -- (136.4366,1006.0333)
      (124.6268,1006.5698) -- (132.5000,1049.6262)
      (172.2209,1031.2390) -- (132.5000,1049.6262)
      (132.5000,919.3839) -- (128.5634,962.9768)
      (136.4366,1006.0333) -- (132.5000,1049.6262);
      \end{tikzpicture}
      \subcaption{reverse double kite}
    \end{minipage}
    \quad
    \begin{minipage}[t]{0.45\linewidth}
      \centering
      \begin{tikzpicture}
      [y=0.5pt, x=0.5pt, yscale=-1.0, xscale=1.0]
      \draw[line width=0.01pt]
      (119.7716,968.4091) -- (97.8865,1006.3153)
      (174.4846,1048.6958) -- (86.9439,1048.6958)
      (108.8290,1010.7896) -- (152.5994,1010.7896)
      (130.7142,1048.6958) -- (108.8290,1010.7896)
      (119.7716,968.4091) -- (108.8290,1010.7896)
      (86.9439,1048.6958) -- (108.8290,1010.7896)
      (130.7142,1048.6958) -- (152.5994,1010.7896)
      (152.5994,1010.7896) -- (174.4846,1048.6958)
      (97.8865,1006.3153) -- (86.9439,1048.6958)
      (24.4813,987.3622) -- (86.9439,1048.6958)
      (97.8865,1006.3153) -- (66.6552,975.6486)
      (55.7126,1018.0290) -- (97.8865,1006.3153)
      (55.7126,1018.0290) -- (66.6552,975.6486)
      (66.6552,975.6486) -- (24.4813,987.3622)
      (88.5404,937.7423) -- (66.6552,975.6486)
      (2.5962,949.4560) -- (46.3665,949.4560)
      (46.3665,949.4560) -- (24.4813,987.3622)
      (88.5404,937.7423) -- (46.3665,949.4560)
      (88.5404,937.7423) -- (119.7716,968.4091)
      (141.6568,968.4091) -- (163.5420,1006.3153)
      (141.6568,968.4091) -- (152.5994,1010.7896)
      (163.5420,1006.3153) -- (174.4846,1048.6958)
      (141.6568,968.4091) -- (152.5994,1010.7896)
      (236.9471,987.3622) -- (174.4846,1048.6958)
      (163.5420,1006.3153) -- (194.7732,975.6485)
      (205.7158,1018.0290) -- (163.5420,1006.3153)
      (205.7158,1018.0290) -- (194.7732,975.6485)
      (194.7732,975.6485) -- (236.9471,987.3622)
      (172.8881,937.7423) -- (194.7732,975.6485)
      (215.0619,949.4560) -- (236.9471,987.3622)
      (172.8881,937.7423) -- (215.0619,949.4560)
      (141.6568,968.4091) -- (172.8881,937.7423)
      (119.7716,968.4091) -- (130.7142,926.0286)
      (130.7142,926.0286) -- (141.6568,968.4091)
      (2.5962,949.4560) -- (24.4813,987.3622)
      (215.0619,949.4560) -- (258.8323,949.4560)
      (258.8323,949.4560) -- (236.9471,987.3622)
      (88.5404,937.7423) -- (130.7142,926.0286)
      (130.7142,926.0286) -- (172.8881,937.7423);
      \end{tikzpicture}
      \subcaption{triplet kite}
    \end{minipage}  
    \caption{Rigid subgraphs}
  \end{figure}
  
  %\noindent
  The triplet kite is a $(2;3;4)$-regular matchstick graph consisting of 22 vertices and 41 edges and has a vertical symmetry. Three triplet kites can be connected together in a way so that the three vertices of the outer triangles become the same. The so formed 4-regular matchstick graph (Fig. 4) offers the possibility to connect each two of the three inner vertices with an additional unit length edge. This property has been used for $n=5$ (Fig. 7) and $n=6$ (Fig. 8). The graphs for $n=9$ (Fig. 11) and $n=11$ (Fig. 13) use additional subgraphs and are the only graphs which contain outside edges, which are not part of an equilateral triangle. The graph for $n=11$, by far the largest graph in this article, has a more complicated geometry, which had to be calculated with a CAS by Stefan Vogel. The graphs for $n=7$ (Fig. 9) are the only flexible graphs and do not consist of the kite-based subgraphs. The geometry of the graphs for $n=5$ with 118 and 120 edges (Fig. 15) combined the subgraphs (a) and (d). The geometry of the graphs for $n=5$ and $n=6$ with 121 edges (Fig. 16 -- 18) is based on the subgraph (d). These graphs are the previous versions of the smallest known follow-up graphs in Figure 7 and 8.
  
  \section{\large{The four smallest known 4-regular matchstick graphs}}
  
  In 1986 Heiko Harborth [2] presented the smallest known example of a 4-regular matchstick graph. This rigid graph has a vertical and a horizontal symmetry.
  
  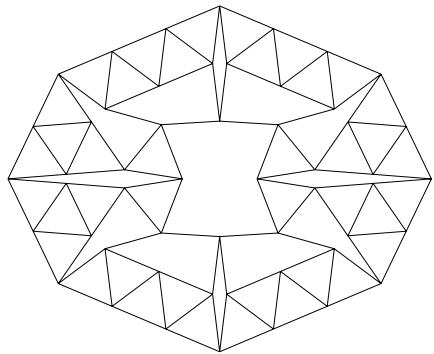
\begin{figure}[!ht]
  \centering
  \begin{tikzpicture}
     [y=0.5pt, x=0.5pt, yscale=-1.0, xscale=1.0]
     \draw[line width=0.01pt]
     (3.6500,918.0914) -- (47.2860,914.6636)
     (47.2860,914.6636) -- (90.9219,911.2359)
     (90.9219,911.2359) -- (66.1354,875.1600)
     (66.1354,875.1600) -- (47.2860,914.6636)
     (47.2860,914.6636) -- (22.4995,878.5877)
     (22.4995,878.5877) -- (3.6500,918.0914)
     (22.4995,878.5877) -- (66.1354,875.1600)
     (66.1354,875.1600) -- (41.3490,839.0840)
     (41.3490,839.0840) -- (22.4995,878.5877)
     (90.9219,911.2359) -- (134.1521,918.0914)
     (134.1521,918.0914) -- (118.4739,877.2252)
     (118.4739,877.2252) -- (90.9219,911.2359)
     (41.3490,839.0840) -- (76.3350,865.3865)
     (76.3350,865.3865) -- (81.6206,821.9364)
     (81.6206,821.9364) -- (41.3490,839.0840)
     (76.3350,865.3865) -- (116.6067,848.2389)
     (116.6067,848.2389) -- (81.6206,821.9364)
     (81.6206,821.9364) -- (121.8923,804.7889)
     (121.8923,804.7889) -- (116.6067,848.2389)
     (116.6067,848.2389) -- (156.8783,831.0913)
     (156.8783,831.0913) -- (121.8923,804.7889)
     (121.8923,804.7889) -- (162.1639,787.6413)
     (162.1639,787.6413) -- (156.8783,831.0913)
     (118.4739,877.2252) -- (76.3350,865.3865)
     (3.6500,918.0914) -- (47.2860,921.5191)
     (47.2860,921.5191) -- (90.9219,924.9468)
     (90.9219,924.9468) -- (66.1354,961.0227)
     (66.1354,961.0227) -- (47.2860,921.5191)
     (47.2860,921.5191) -- (22.4995,957.5950)
     (22.4995,957.5950) -- (3.6500,918.0914)
     (22.4995,957.5950) -- (66.1354,961.0228)
     (66.1354,961.0227) -- (41.3490,997.0987)
     (41.3490,997.0987) -- (22.4995,957.5950)
     (90.9219,924.9468) -- (134.1521,918.0914)
     (134.1521,918.0914) -- (118.4739,958.9575)
     (118.4739,958.9575) -- (90.9219,924.9468)
     (41.3490,997.0987) -- (76.3350,970.7962)
     (76.3350,970.7962) -- (81.6206,1014.2463)
     (81.6206,1014.2463) -- (41.3490,997.0987)
     (76.3350,970.7962) -- (116.6067,987.9438)
     (116.6067,987.9438) -- (81.6206,1014.2463)
     (81.6206,1014.2463) -- (121.8923,1031.3939)
     (121.8923,1031.3939) -- (116.6067,987.9438)
     (116.6067,987.9438) -- (156.8783,1005.0914)
     (156.8783,1005.0914) -- (121.8922,1031.3939)
     (121.8923,1031.3939) -- (162.1639,1048.5415)
     (162.1639,1048.5415) -- (156.8783,1005.0914)
     (118.4739,958.9575) -- (76.3350,970.7962)
     (320.6778,918.0914) -- (277.0418,914.6636)
     (277.0418,914.6636) -- (233.4059,911.2359)
     (233.4059,911.2359) -- (258.1923,875.1600)
     (258.1923,875.1600) -- (277.0418,914.6637)
     (277.0418,914.6636) -- (301.8283,878.5877)
     (301.8283,878.5877) -- (320.6778,918.0914)
     (301.8283,878.5877) -- (258.1923,875.1600)
     (258.1923,875.1600) -- (282.9788,839.0840)
     (282.9788,839.0840) -- (301.8283,878.5877)
     (233.4059,911.2359) -- (190.1757,918.0914)
     (190.1757,918.0914) -- (205.8538,877.2252)
     (205.8538,877.2252) -- (233.4059,911.2359)
     (282.9788,839.0840) -- (247.9927,865.3865)
     (247.9927,865.3865) -- (242.7071,821.9364)
     (242.7072,821.9364) -- (282.9788,839.0840)
     (247.9927,865.3865) -- (207.7211,848.2389)
     (207.7211,848.2389) -- (242.7072,821.9364)
     (242.7072,821.9364) -- (202.4355,804.7889)
     (202.4355,804.7889) -- (207.7211,848.2389)
     (207.7211,848.2389) -- (167.4495,831.0913)
     (167.4495,831.0913) -- (202.4355,804.7889)
     (202.4355,804.7889) -- (162.1639,787.6413)
     (162.1639,787.6413) -- (167.4495,831.0913)
     (205.8538,877.2252) -- (247.9927,865.3865)
     (320.6778,918.0914) -- (277.0418,921.5191)
     (277.0418,921.5191) -- (233.4059,924.9468)
     (233.4059,924.9468) -- (258.1923,961.0227)
     (258.1923,961.0227) -- (277.0418,921.5191)
     (277.0418,921.5191) -- (301.8283,957.5950)
     (301.8283,957.5950) -- (320.6778,918.0914)
     (301.8283,957.5950) -- (258.1923,961.0228)
     (258.1923,961.0227) -- (282.9788,997.0987)
     (282.9788,997.0987) -- (301.8283,957.5950)
     (233.4059,924.9468) -- (190.1757,918.0914)
     (190.1757,918.0914) -- (205.8538,958.9575)
     (205.8538,958.9575) -- (233.4059,924.9468)
     (282.9788,997.0987) -- (247.9927,970.7962)
     (247.9927,970.7962) -- (242.7071,1014.2463)
     (242.7072,1014.2463) -- (282.9788,997.0987)
     (247.9927,970.7962) -- (207.7211,987.9438)
     (207.7211,987.9438) -- (242.7072,1014.2463)
     (242.7072,1014.2463) -- (202.4355,1031.3939)
     (202.4355,1031.3939) -- (207.7211,987.9438)
     (207.7211,987.9438) -- (167.4495,1005.0914)
     (167.4495,1005.0914) -- (202.4355,1031.3939)
     (202.4355,1031.3939) -- (162.1639,1048.5415)
     (162.1639,1048.5415) -- (167.4495,1005.0914)
     (205.8538,958.9575) -- (247.9927,970.7962)
     (156.8783,831.0913) -- (162.1639,874.5414)
     (162.1639,874.5414) -- (167.4495,831.0913)
     (156.8783,1005.0914) -- (162.1639,961.6413)
     (162.1639,961.6413) -- (167.4495,1005.0914)
     (118.4739,958.9575) -- (162.1639,961.6413)
     (162.1639,961.6413) -- (205.8538,958.9575)
     (118.4739,877.2252) -- (162.1639,874.5414)
     (162.1639,874.5414) -- (205.8538,877.2252);
     \end{tikzpicture}
     \caption{4-regular matchstick graph with 52 vertices and 104 edges.}
  \end{figure}
  
  \newpage
  
  \noindent On July 3, 2016 the authors discovered a new second smallest known example of a 4-regular matchstick graph. This rigid graph has a vertical and a horizontal symmetry and is based on the Harborth graph. The history of this graph is a little bit intricate and begins on April 24, 2016 [8].
  
  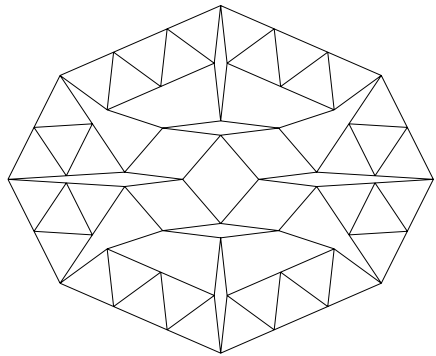
\begin{figure}[!ht]
  	\centering
    \begin{tikzpicture}
     [y=0.5pt, x=0.5pt, yscale=-1.0, xscale=1.0]
     \draw[line width=0.01pt]
     (46.5102,993.8452) -- (81.7775,967.9211)
     (81.7775,967.9211) -- (86.5948,1011.4255)
     (86.5948,1011.4255) -- (46.5102,993.8452)
     (81.7775,967.9211) -- (121.8621,985.5015)
     (121.8621,985.5014) -- (86.5948,1011.4255)
     (86.5948,1011.4255) -- (126.6794,1029.0059)
     (126.6794,1029.0059) -- (121.8621,985.5014)
     (121.8621,985.5014) -- (161.9467,1003.0818)
     (161.9467,1003.0818) -- (126.6794,1029.0059)
     (126.6794,1029.0059) -- (166.7639,1046.5862)
     (166.7639,1046.5862) -- (161.9467,1003.0818)
     (46.5102,993.8452) -- (70.7013,957.3674)
     (70.7014,957.3674) -- (94.8926,920.8896)
     (27.0151,954.6561) -- (70.7013,957.3674)
     (70.7014,957.3674) -- (51.2062,918.1783)
     (51.2062,918.1783) -- (27.0151,954.6561)
     (46.5102,993.8452) -- (27.0151,954.6561)
     (27.0151,954.6561) -- (7.5199,915.4670)
     (94.8926,920.8896) -- (51.2062,918.1783)
     (51.2062,918.1783) -- (7.5199,915.4670)
     (81.7775,967.9211) -- (123.3297,954.1638)
     (123.3297,954.1638) -- (94.8926,920.8896)
     (161.9467,1003.0818) -- (166.7639,959.5774)
     (166.7639,959.5774) -- (123.3297,954.1638)
     (287.0177,993.8452) -- (251.7504,967.9211)
     (251.7504,967.9211) -- (246.9331,1011.4255)
     (246.9331,1011.4255) -- (287.0177,993.8452)
     (251.7504,967.9211) -- (211.6658,985.5015)
     (211.6658,985.5014) -- (246.9331,1011.4255)
     (246.9331,1011.4255) -- (206.8485,1029.0059)
     (206.8485,1029.0059) -- (211.6658,985.5015)
     (211.6658,985.5014) -- (171.5812,1003.0818)
     (171.5812,1003.0818) -- (206.8485,1029.0059)
     (206.8485,1029.0059) -- (166.7640,1046.5862)
     (166.7639,1046.5862) -- (171.5812,1003.0818)
     (287.0177,993.8452) -- (262.8265,957.3674)
     (262.8265,957.3674) -- (238.6353,920.8896)
     (306.5128,954.6561) -- (262.8265,957.3674)
     (262.8265,957.3674) -- (282.3216,918.1783)
     (282.3216,918.1783) -- (306.5128,954.6561)
     (287.0177,993.8452) -- (306.5128,954.6561)
     (306.5128,954.6561) -- (326.0079,915.4670)
     (238.6353,920.8896) -- (282.3216,918.1783)
     (282.3216,918.1783) -- (326.0079,915.4670)
     (251.7504,967.9211) -- (210.1982,954.1638)
     (210.1982,954.1638) -- (238.6353,920.8896)
     (171.5812,1003.0818) -- (166.7639,959.5774)
     (166.7639,959.5774) -- (210.1982,954.1638)
     (123.3297,954.1638) -- (166.7639,948.7503)
     (166.7639,948.7503) -- (210.1982,954.1638)
     (46.5102,837.0888) -- (81.7775,863.0129)
     (81.7775,863.0129) -- (86.5948,819.5085)
     (86.5948,819.5085) -- (46.5102,837.0888)
     (81.7775,863.0129) -- (121.8621,845.4326)
     (121.8621,845.4326) -- (86.5948,819.5085)
     (86.5948,819.5085) -- (126.6794,801.9281)
     (126.6794,801.9281) -- (121.8621,845.4326)
     (121.8621,845.4326) -- (161.9467,827.8522)
     (161.9467,827.8522) -- (126.6794,801.9281)
     (126.6794,801.9281) -- (166.7639,784.3477)
     (166.7639,784.3477) -- (161.9467,827.8522)
     (46.5102,837.0888) -- (70.7013,873.5666)
     (70.7014,873.5666) -- (94.8926,910.0444)
     (27.0151,876.2779) -- (70.7013,873.5666)
     (70.7014,873.5666) -- (51.2063,912.7557)
     (51.2062,912.7557) -- (27.0151,876.2779)
     (46.5102,837.0888) -- (27.0151,876.2779)
     (27.0151,876.2779) -- (7.5199,915.4670)
     (94.8926,910.0444) -- (51.2062,912.7557)
     (51.2062,912.7557) -- (7.5199,915.4670)
     (81.7775,863.0129) -- (123.3297,876.7702)
     (123.3297,876.7702) -- (94.8926,910.0444)
     (161.9467,827.8522) -- (166.7639,871.3567)
     (166.7639,871.3566) -- (123.3297,876.7702)
     (287.0177,837.0888) -- (251.7504,863.0129)
     (251.7504,863.0129) -- (246.9331,819.5085)
     (246.9331,819.5085) -- (287.0177,837.0888)
     (251.7504,863.0129) -- (211.6658,845.4326)
     (211.6658,845.4326) -- (246.9331,819.5085)
     (246.9331,819.5085) -- (206.8485,801.9281)
     (206.8485,801.9281) -- (211.6658,845.4326)
     (211.6658,845.4326) -- (171.5812,827.8522)
     (171.5812,827.8522) -- (206.8485,801.9281)
     (206.8485,801.9281) -- (166.7639,784.3477)
     (166.7639,784.3477) -- (171.5812,827.8522)
     (287.0177,837.0888) -- (262.8265,873.5666)
     (262.8265,873.5666) -- (238.6353,910.0444)
     (306.5128,876.2779) -- (262.8265,873.5666)
     (262.8265,873.5666) -- (282.3216,912.7557)
     (282.3216,912.7557) -- (306.5128,876.2779)
     (287.0177,837.0888) -- (306.5128,876.2779)
     (306.5128,876.2779) -- (326.0079,915.4670)
     (238.6353,910.0444) -- (282.3216,912.7557)
     (282.3216,912.7557) -- (326.0079,915.4670)
     (251.7504,863.0129) -- (210.1982,876.7702)
     (210.1982,876.7702) -- (238.6353,910.0444)
     (171.5812,827.8522) -- (166.7639,871.3567)
     (166.7639,871.3566) -- (210.1982,876.7702)
     (123.3297,876.7702) -- (166.7639,882.1837)
     (166.7639,882.1837) -- (210.1982,876.7702)
     (94.8926,920.8896) -- (138.3257,915.4670)
     (138.3257,915.4670) -- (94.8926,910.0444)
     (238.6353,920.8896) -- (195.2022,915.4670)
     (195.2022,915.4670) -- (238.6353,910.0444)
     (138.3257,915.4670) -- (166.7639,948.7503)
     (166.7639,948.7503) -- (195.2022,915.4670)
     (195.2022,915.4670) -- (166.7639,882.1837)
     (166.7639,882.1837) -- (138.3257,915.4670);
    \end{tikzpicture}
    \caption{4-regular matchstick graph with 54 vertices and 108 edges.}
  \end{figure}
  
  \noindent\\ On April 15, 2016 Mike Winkler discovered a new third smallest known example of a 4-regular matchstick graph [9][10]. This rigid graph has a rotational symmetry of order 3 and consists of three overlapping triplet kites.
  
  \begin{figure}[!ht]
  \centering
  \begin{tikzpicture}
     [y=0.5pt, x=0.5pt, yscale=-1.0, xscale=1.0]
     \draw[line width=0.01pt]
     (146.4577,967.9964) -- (124.5725,1005.9026)
     (201.1706,1048.2831) -- (113.6299,1048.2831)
     (135.5151,1010.3768) -- (179.2855,1010.3768)
     (157.4003,1048.2831) -- (135.5151,1010.3768)
     (146.4577,967.9964) -- (135.5151,1010.3768)
     (113.6299,1048.2831) -- (135.5151,1010.3768)
     (157.4003,1048.2831) -- (179.2855,1010.3768)
     (179.2855,1010.3768) -- (201.1706,1048.2831)
     (124.5725,1005.9026) -- (113.6300,1048.2831)
     (51.1674,986.9495) -- (113.6300,1048.2831)
     (124.5725,1005.9026) -- (93.3413,975.2358)
     (82.3987,1017.6163) -- (124.5725,1005.9026)
     (82.3987,1017.6163) -- (93.3413,975.2358)
     (93.3413,975.2358) -- (51.1674,986.9495)
     (115.2264,937.3296) -- (93.3413,975.2358)
     (29.2822,949.0432) -- (73.0526,949.0432)
     (73.0526,949.0432) -- (51.1674,986.9495)
     (115.2264,937.3296) -- (73.0526,949.0432)
     (115.2264,937.3296) -- (146.4577,967.9964)
     (168.3429,967.9964) -- (190.2281,1005.9026)
     (179.2855,1010.3768) -- (135.5151,1010.3768)
     (168.3429,967.9964) -- (179.2855,1010.3768)
     (190.2281,1005.9026) -- (201.1706,1048.2831)
     (263.6332,986.9495) -- (201.1707,1048.2831)
     (190.2281,1005.9026) -- (221.4593,975.2358)
     (232.4019,1017.6163) -- (190.2281,1005.9026)
     (232.4019,1017.6163) -- (221.4593,975.2358)
     (221.4593,975.2358) -- (263.6332,986.9495)
     (199.5741,937.3296) -- (221.4593,975.2358)
     (241.7480,949.0432) -- (263.6332,986.9495)
     (199.5741,937.3296) -- (241.7480,949.0432)
     (168.3429,967.9964) -- (199.5741,937.3296)
     (146.4577,967.9964) -- (157.4003,925.6159)
     (157.4003,925.6159) -- (168.3429,967.9964)
     (115.2264,937.3296) -- (157.4003,925.6159)
     (29.2822,949.0432) -- (51.1674,986.9495)
     (157.4003,925.6159) -- (199.5741,937.3296)
     (241.7480,949.0432) -- (285.5183,949.0432)
     (285.5183,949.0432) -- (263.6332,986.9495)
     (210.5167,857.0429) -- (232.4019,819.1366)
     (307.4035,864.2823) -- (263.6332,788.4698)
     (241.7480,826.3761) -- (263.6332,864.2823)
     (285.5183,826.3761) -- (241.7480,826.3761)
     (210.5167,857.0429) -- (241.7480,826.3761)
     (263.6332,788.4698) -- (241.7480,826.3761)
     (285.5183,826.3761) -- (263.6332,864.2823)
     (263.6332,864.2823) -- (307.4035,864.2823)
     (232.4019,819.1366) -- (263.6332,788.4699)
     (179.2855,765.0425) -- (263.6332,788.4699)
     (232.4019,819.1367) -- (190.2281,807.4230)
     (221.4593,776.7562) -- (232.4019,819.1366)
     (221.4593,776.7562) -- (190.2281,807.4230)
     (190.2281,807.4230) -- (179.2855,765.0425)
     (168.3429,845.3292) -- (190.2281,807.4230)
     (135.5151,765.0425) -- (157.4003,802.9488)
     (157.4003,802.9488) -- (179.2855,765.0425)
     (168.3429,845.3292) -- (157.4003,802.9488)
     (168.3429,845.3292) -- (210.5167,857.0429)
     (221.4593,875.9960) -- (265.2297,875.9960)
     (221.4593,875.9960) -- (263.6332,864.2823)
     (307.4035,864.2823) -- (263.6332,864.2823)
     (265.2297,875.9960) -- (307.4035,864.2823)
     (221.4593,875.9960) -- (263.6332,864.2823)
     (285.5183,949.0432) -- (307.4035,864.2823)
     (265.2297,875.9960) -- (254.2871,918.3764)
     (296.4609,906.6628) -- (265.2297,875.9960)
     (254.2871,918.3764) -- (285.5183,949.0432)
     (210.5167,918.3764) -- (254.2871,918.3764)
     (210.5167,918.3764) -- (241.7480,949.0432)
     (221.4593,875.9960) -- (210.5167,918.3764)
     (210.5167,857.0429) -- (179.2855,887.7097)
     (179.2855,887.7097) -- (221.4593,875.9960)
     (168.3429,845.3292) -- (179.2855,887.7097)
     (135.5151,765.0425) -- (179.2855,765.0425)
     (179.2855,887.7097) -- (210.5167,918.3764)
     (93.3413,875.9960) -- (49.5709,875.9960)
     (51.1674,788.4699) -- (7.3970,864.2823)
     (51.1674,864.2823) -- (73.0526,826.3761)
     (29.2822,826.3761) -- (51.1674,864.2823)
     (93.3413,875.9960) -- (51.1674,864.2823)
     (7.3970,864.2823) -- (51.1674,864.2823)
     (29.2822,826.3761) -- (73.0526,826.3761)
     (73.0526,826.3761) -- (51.1674,788.4699)
     (49.5709,875.9960) -- (7.3970,864.2823)
     (29.2822,949.0432) -- (7.3971,864.2823)
     (49.5709,875.9960) -- (60.5135,918.3765)
     (18.3396,906.6628) -- (49.5709,875.9960)
     (18.3396,906.6628) -- (60.5135,918.3765)
     (60.5135,918.3765) -- (29.2822,949.0432)
     (104.2839,918.3765) -- (60.5135,918.3765)
     (73.0526,949.0432) -- (29.2822,949.0432)
     (104.2839,918.3765) -- (73.0526,949.0432)
     (104.2839,918.3765) -- (93.3413,875.9960)
     (104.2838,857.0429) -- (82.3987,819.1366)
     (104.2838,857.0429) -- (73.0526,826.3761)
     (51.1674,788.4699) -- (73.0526,826.3761)
     (82.3987,819.1366) -- (51.1674,788.4699)
     (104.2838,857.0429) -- (73.0526,826.3761)
     (135.5151,765.0425) -- (51.1674,788.4699)
     (82.3987,819.1366) -- (124.5725,807.4230)
     (93.3412,776.7562) -- (82.3987,819.1366)
     (93.3412,776.7562) -- (124.5725,807.4230)
     (124.5725,807.4230) -- (135.5151,765.0425)
     (146.4577,845.3292) -- (124.5725,807.4230)
     (146.4577,845.3292) -- (157.4003,802.9487)
     (104.2838,857.0429) -- (146.4577,845.3292)
     (146.4577,845.3292) -- (135.5151,887.7097)
     (135.5151,887.7097) -- (104.2838,857.0429)
     (93.3413,875.9960) -- (135.5151,887.7097)
     (135.5151,887.7097) -- (104.2839,918.3765)
     (254.2871,918.3764) -- (296.4609,906.6628);
    \end{tikzpicture}
    \caption{4-regular matchstick graph with 57 vertices and 114 edges.}
  \end{figure}
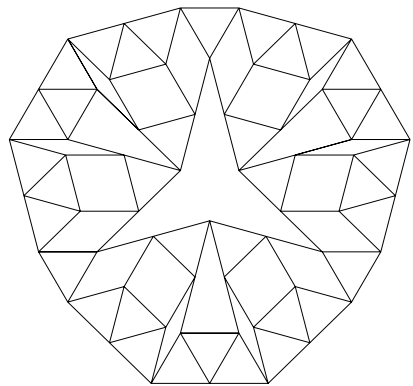
  
  \noindent The graph in Figure 4 is also the basis for the smallest known $(4;n)$-regular matchstick graphs for $n=5$ and $n=6$.
  
  \newpage
  
  \noindent From 1986 to April 15, 2016 the second smallest known example of a 4-regular matchstick graph consisted of 60 vertices and 120 edges.
  
  \begin{figure}[!ht]
    \centering
    \begin{minipage}[t]{0.45\linewidth}
      \centering
      \begin{tikzpicture}
       [y=0.47pt, x=0.47pt, yscale=-1.0, xscale=1.0]
       \draw[line width=0.01pt,red] % red edges left
       (220.3499,1038.5634) -- (132.8092,1038.5634)
       (154.6944,1000.6571) -- (198.4647,1000.6571)
       (176.5796,1038.5634) -- (154.6944,1000.6571)
       (132.8092,1038.5634) -- (154.6944,1000.6571)
       (176.5796,1038.5634) -- (198.4647,1000.6571)
       (198.4647,1000.6571) -- (220.3499,1038.5634)
       (198.4647,1000.6571) -- (154.6944,1000.6571)
       (154.6944,1000.6571) -- (154.6944,956.8868)
       (154.6944,956.8868) -- (198.4647,956.8868)
       (198.4647,956.8868) -- (198.4647,1000.6571);
       \draw[line width=0.01pt,black] % black edges left
       (220.3499,1038.5634) -- (296.1624,994.7930)
       (258.2561,972.9078) -- (220.3499,994.7930)
       (258.2561,1016.6782) -- (258.2561,972.9078)
       (296.1624,994.7930) -- (258.2561,972.9078)
       (258.2561,1016.6782) -- (220.3499,994.7930)
       (220.3499,994.7930) -- (220.3499,1038.5634)
       (220.3499,994.7930) -- (258.2561,972.9078)
       (258.2561,972.9078) -- (236.3709,935.0016)
       (236.3709,935.0016) -- (198.4647,956.8868)
       (198.4647,956.8868) -- (220.3499,994.7930)
       (56.9967,994.7930) -- (13.2264,918.9805)
       (56.9967,918.9805) -- (78.8819,956.8868)
       (35.1115,956.8868) -- (56.9967,918.9805)
       (13.2264,918.9805) -- (56.9967,918.9805)
       (35.1115,956.8868) -- (78.8819,956.8868)
       (78.8819,956.8868) -- (56.9967,994.7930)
       (78.8819,956.8868) -- (56.9967,918.9805)
       (56.9967,918.9805) -- (94.9029,897.0954)
       (94.9029,897.0954) -- (116.7881,935.0016)
       (116.7881,935.0016) -- (78.8819,956.8868)
       (56.9967,994.7930) -- (132.8092,1038.5634)
       (132.8092,994.7930) -- (94.9030,972.9078)
       (94.9030,1016.6782) -- (132.8092,994.7930)
       (132.8092,1038.5634) -- (132.8092,994.7930)
       (94.9030,1016.6782) -- (94.9030,972.9078)
       (94.9030,972.9078) -- (56.9967,994.7930)
       (94.9030,972.9078) -- (132.8092,994.7930)
       (132.8092,994.7930) -- (154.6944,956.8868)
       (154.6944,956.8868) -- (116.7881,935.0016)
       (116.7881,935.0016) -- (94.9030,972.9078)
       (13.2264,831.4398) -- (56.9967,755.6274)
       (78.8819,793.5336) -- (56.9967,831.4398)
       (35.1115,793.5336) -- (78.8819,793.5336)
       (56.9967,755.6274) -- (78.8819,793.5336)
       (35.1115,793.5336) -- (56.9967,831.4398)
       (56.9967,831.4398) -- (13.2264,831.4398)
       (56.9967,831.4398) -- (78.8819,793.5336)
       (78.8819,793.5336) -- (116.7881,815.4188)
       (116.7881,815.4188) -- (94.9029,853.3250)
       (94.9029,853.3250) -- (56.9967,831.4398)
       (13.2264,831.4398) -- (13.2264,918.9805)
       (51.1326,897.0954) -- (51.1326,853.3250)
       (13.2264,875.2102) -- (51.1326,897.0954)
       (13.2264,918.9805) -- (51.1326,897.0954)
       (13.2264,875.2102) -- (51.1326,853.3250)
       (51.1326,853.3250) -- (13.2264,831.4398)
       (51.1326,853.3250) -- (51.1326,897.0954)
       (51.1326,897.0954) -- (94.9029,897.0954)
       (94.9029,897.0954) -- (94.9029,853.3250)
       (94.9029,853.3250) -- (51.1326,853.3250)
       (132.8092,711.8570) -- (220.3499,711.8570)
       (198.4647,749.7632) -- (154.6944,749.7632)
       (176.5795,711.8570) -- (198.4647,749.7632)
       (220.3499,711.8570) -- (198.4647,749.7632)
       (176.5795,711.8570) -- (154.6944,749.7632)
       (154.6944,749.7632) -- (132.8092,711.8570)
       (154.6944,749.7632) -- (198.4647,749.7632)
       (198.4647,749.7632) -- (198.4647,793.5336)
       (198.4647,793.5336) -- (154.6944,793.5336)
       (154.6944,793.5336) -- (154.6944,749.7632)
       (132.8092,711.8570) -- (56.9967,755.6274)
       (94.9029,777.5125) -- (132.8092,755.6274)
       (94.9029,733.7422) -- (94.9029,777.5125)
       (56.9967,755.6274) -- (94.9029,777.5125)
       (94.9029,733.7422) -- (132.8092,755.6274)
       (132.8092,755.6274) -- (132.8092,711.8570)
       (132.8092,755.6274) -- (94.9029,777.5126)
       (94.9029,777.5125) -- (116.7881,815.4188)
       (116.7881,815.4188) -- (154.6944,793.5336)
       (154.6944,793.5336) -- (132.8092,755.6274)
       (296.1624,755.6274) -- (339.9327,831.4398)
       (296.1624,831.4398) -- (274.2772,793.5336)
       (318.0475,793.5336) -- (296.1624,831.4398)
       (339.9327,831.4398) -- (296.1624,831.4398)
       (318.0475,793.5336) -- (274.2772,793.5336)
       (274.2772,793.5336) -- (296.1624,755.6274)
       (274.2772,793.5336) -- (296.1624,831.4398)
       (296.1624,831.4398) -- (258.2561,853.3250)
       (258.2561,853.3250) -- (236.3709,815.4188)
       (236.3709,815.4188) -- (274.2772,793.5336)
       (296.1624,755.6274) -- (220.3499,711.8570)
       (220.3499,755.6274) -- (258.2561,777.5125)
       (258.2561,733.7422) -- (220.3499,755.6274)
       (220.3499,711.8570) -- (220.3499,755.6274)
       (258.2561,733.7422) -- (258.2561,777.5125)
       (258.2561,777.5125) -- (296.1624,755.6274)
       (258.2561,777.5125) -- (220.3499,755.6274)
       (220.3499,755.6274) -- (198.4647,793.5336)
       (198.4647,793.5336) -- (236.3709,815.4188)
       (236.3709,815.4188) -- (258.2561,777.5125)
       (339.9327,918.9805) -- (296.1624,994.7930)
       (274.2772,956.8868) -- (296.1624,918.9805)
       (318.0475,956.8868) -- (274.2772,956.8868)
       (296.1624,994.7930) -- (274.2772,956.8868)
       (318.0475,956.8868) -- (296.1624,918.9805)
       (296.1624,918.9805) -- (339.9327,918.9805)
       (296.1624,918.9805) -- (274.2772,956.8868)
       (274.2772,956.8868) -- (236.3710,935.0016)
       (236.3710,935.0016) -- (258.2561,897.0953)
       (258.2561,897.0954) -- (296.1624,918.9805)
       (339.9327,918.9805) -- (339.9327,831.4398)
       (302.0265,853.3250) -- (302.0265,897.0954)
       (339.9327,875.2102) -- (302.0265,853.3250)
       (339.9327,831.4398) -- (302.0265,853.3250)
       (339.9327,875.2102) -- (302.0265,897.0953)
       (302.0265,897.0954) -- (339.9327,918.9805)
       (302.0265,897.0954) -- (302.0265,853.3250)
       (302.0265,853.3250) -- (258.2561,853.3250)
       (258.2561,853.3250) -- (258.2561,897.0954)
       (258.2561,897.0954) -- (302.0265,897.0954);
      \end{tikzpicture}
      \subcaption*{v1}
    \end{minipage}
    \quad
    \begin{minipage}[t]{0.45\linewidth}
      \centering
      \begin{tikzpicture}
       [y=0.47pt, x=0.47pt, yscale=-1.0, xscale=1.0]
       \draw[line width=0.01pt,red] % red edges right
       (524.9146,958.2952) -- (503.0294,996.2015)
       (579.6276,1038.5819) -- (492.0869,1038.5819)
       (513.9720,1000.6757) -- (557.7424,1000.6757)
       (535.8572,1038.5819) -- (513.9720,1000.6757)
       (524.9146,958.2952) -- (513.9720,1000.6757)
       (492.0869,1038.5819) -- (513.9720,1000.6757)
       (535.8572,1038.5819) -- (557.7424,1000.6757)
       (557.7424,1000.6757) -- (579.6276,1038.5819)
       (503.0294,996.2015) -- (492.0869,1038.5819)
       (429.6243,977.2483) -- (492.0869,1038.5819)
       (503.0294,996.2014) -- (471.7982,965.5347)
       (460.8556,1007.9151) -- (503.0294,996.2014)
       (460.8556,1007.9151) -- (471.7982,965.5347)
       (471.7982,965.5347) -- (429.6243,977.2483)
       (493.6834,927.6284) -- (471.7982,965.5347)
       (493.6834,927.6284) -- (524.9146,958.2952)
        (524.9146,958.2952) -- (568.6850,958.2952)
       (568.6850,958.2952) -- (557.7424,1000.6757);      
       \draw[line width=0.01pt,black] % black edges right
       (493.6834,927.6284) -- (451.5095,939.3421)
       (471.7982,889.7222) -- (428.0278,889.7222)
       (429.6243,977.2483) -- (385.8540,901.4358)
       (429.6243,901.4358) -- (451.5095,939.3421)
       (407.7392,939.3421) -- (429.6243,901.4358)
       (471.7982,889.7222) -- (429.6243,901.4358)
       (385.8540,901.4358) -- (429.6243,901.4358)
       (407.7392,939.3421) -- (451.5095,939.3421)
       (451.5095,939.3421) -- (429.6243,977.2483)
       (428.0278,889.7222) -- (385.8540,901.4358)
       (407.7392,816.6749) -- (385.8540,901.4359)
       (428.0278,889.7222) -- (438.9704,847.3417)
       (396.7966,859.0554) -- (428.0278,889.7222)
       (396.7966,859.0554) -- (438.9704,847.3417)
       (438.9704,847.3417) -- (407.7392,816.6749)
       (482.7408,847.3417) -- (438.9704,847.3417)
       (482.7408,847.3417) -- (451.5095,816.6749)
       (482.7408,847.3417) -- (471.7982,889.7222)
       (451.5095,939.3421) -- (429.6243,901.4358)
       (504.6259,809.4355) -- (482.7408,771.5293)
       (407.7392,816.6749) -- (451.5095,740.8625)
       (473.3947,778.7687) -- (451.5095,816.6749)
       (429.6243,778.7687) -- (473.3947,778.7687)
       (504.6259,809.4355) -- (473.3947,778.7687)
       (451.5095,740.8625) -- (473.3947,778.7687)
       (429.6243,778.7687) -- (451.5095,816.6749)
       (451.5095,816.6749) -- (407.7392,816.6749)
       (482.7408,771.5293) -- (451.5095,740.8625)
       (535.8572,717.4351) -- (451.5095,740.8625)
       (482.7408,771.5293) -- (524.9146,759.8156)
       (493.6834,729.1488) -- (482.7408,771.5293)
       (493.6834,729.1488) -- (524.9146,759.8156)
       (524.9146,759.8156) -- (535.8572,717.4351)
       (546.7998,797.7218) -- (524.9146,759.8156)
       (546.7998,797.7218) -- (557.7424,755.3414)
       (546.7998,797.7218) -- (504.6259,809.4355)
       (451.5095,816.6749) -- (473.3947,778.7687)
       (590.5701,797.7218) -- (612.4553,759.8156)
       (535.8572,717.4351) -- (623.3979,717.4351)
       (601.5127,755.3414) -- (557.7424,755.3414)
       (579.6276,717.4351) -- (601.5127,755.3414)
       (590.5701,797.7218) -- (601.5127,755.3414)
       (623.3979,717.4351) -- (601.5127,755.3414)
       (579.6276,717.4351) -- (557.7424,755.3414)
       (557.7424,755.3414) -- (535.8572,717.4351)
       (612.4553,759.8156) -- (623.3979,717.4351)
       (685.8604,778.7687) -- (623.3979,717.4351)
       (612.4553,759.8156) -- (643.6866,790.4824)
       (654.6292,748.1019) -- (612.4553,759.8156)
       (654.6292,748.1019) -- (643.6866,790.4824)
       (643.6866,790.4824) -- (685.8604,778.7687)
       (621.8014,828.3886) -- (643.6866,790.4824)
       (621.8014,828.3886) -- (663.9752,816.6749)
       (621.8014,828.3886) -- (590.5701,797.7218)
       (557.7424,755.3414) -- (601.5127,755.3414)
       (643.6866,866.2948) -- (687.4569,866.2948)
       (685.8604,778.7687) -- (729.6308,854.5811)
       (685.8604,854.5812) -- (663.9752,816.6749)
       (707.7456,816.6749) -- (685.8604,854.5812)
       (643.6866,866.2948) -- (685.8604,854.5812)
       (729.6308,854.5811) -- (685.8604,854.5812)
       (707.7456,816.6749) -- (663.9752,816.6749)
       (663.9752,816.6749) -- (685.8604,778.7687)
       (687.4569,866.2948) -- (729.6308,854.5812)
       (707.7456,939.3421) -- (729.6308,854.5811)
       (687.4569,866.2948) -- (676.5143,908.6753)
       (718.6882,896.9616) -- (687.4569,866.2948)
       (718.6882,896.9616) -- (676.5143,908.6753)
       (676.5143,908.6753) -- (707.7456,939.3421)
       (632.7440,908.6753) -- (676.5143,908.6753)
       (632.7440,908.6753) -- (663.9753,939.3421)
       (632.7440,908.6753) -- (643.6866,866.2948)
       (663.9752,816.6749) -- (685.8604,854.5812)
       (610.8588,946.5815) -- (632.7440,984.4878)
       (707.7456,939.3421) -- (663.9753,1015.1546)
       (642.0901,977.2483) -- (663.9752,939.3421)
       (685.8604,977.2483) -- (642.0901,977.2483)
       (610.8588,946.5815) -- (642.0901,977.2483)
       (663.9753,1015.1546) -- (642.0901,977.2483)
       (685.8604,977.2483) -- (663.9752,939.3421)
       (663.9752,939.3421) -- (707.7456,939.3421)
       (632.7440,984.4878) -- (663.9753,1015.1546)
       (579.6276,1038.5819) -- (663.9753,1015.1545)
       (632.7440,984.4878) -- (590.5701,996.2014)
       (621.8014,1026.8682) -- (632.7440,984.4878)
       (621.8014,1026.8682) -- (590.5701,996.2014)
       (590.5701,996.2014) -- (579.6276,1038.5819)
       (568.6850,958.2952) -- (590.5701,996.2014)
       (568.6850,958.2952) -- (610.8588,946.5815)
       (663.9752,939.3421) -- (642.0901,977.2483)
       (471.7982,889.7222) -- (493.6834,927.6284)
       (482.7408,847.3417) -- (504.6259,809.4355)
       (546.7998,797.7218) -- (590.5701,797.7218)
       (621.8014,828.3886) -- (643.6866,866.2948)
       (610.8588,946.5815) -- (632.7440,908.6753);
       % big nodes right
       \fill(579.6276,1038.5819) circle (1.5pt);
       \fill(503.0294,996.2015) circle (1.5pt);  
       % character a and b right
       \coordinate[label=right:$a$] (a) at (580,1025);
       \coordinate[label=left:$b$] (b) at (514,975);
      \end{tikzpicture}
      \subcaption*{v2}
    \end{minipage}
  \caption{4-regular matchstick graphs with 60 vertices and 120 edges.}
  \end{figure}
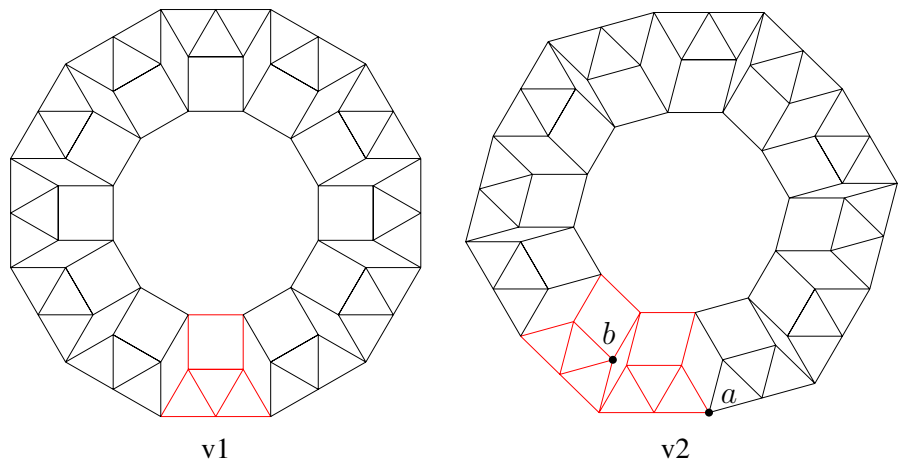
  
  \noindent The graphs in Figure 5 are flexible and each of them can be transformed into the other. Whereby the graph v2 shows only one possibility of transforming. The graph v1 has a rotational symmetry of order 12 and consists of 12 identical subgraphs consisting of 7 vertices and 10 edges (red). Each transformed version, as the graph v2, has a rotational symmetry of order 6 and consists of 6 identical subgraphs consisting of 12 vertices and 20 edges (red).
  
  The transformation in the graph v2 has been chosen so, that the line segment between the vertices $a$ and $b$ measures exactly two unit lengths. This kind of transforming is exactly the one we need to construct the triplet kite with (Fig. 6). This shows the close geometric relationship between the graphs in Figure 4 and Figure 5.
  
  \begin{figure}[!ht]
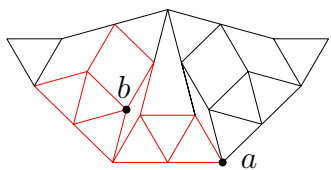

    \centering

    \caption{The geometry of the triplet kite.}
  \end{figure}
  
  \newpage
  
  \section{\large{The smallest known $(4;n)$-regular matchstick graphs for $5\leq n\leq11$}}
  
  \begin{figure}[!ht]
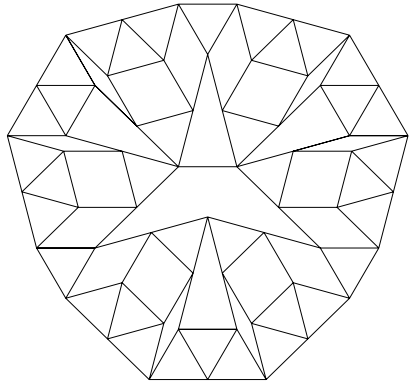

  	\centering
    \caption{$(4;5)$-regular matchstick graph with 57 vertices and 115 edges.}
  	%
  \end{figure}
  \noindent Discovered April 15, 2016 by M. Winkler. This rigid graph has a vertical symmetry and contains three overlapped triplet kites.
  
  \begin{center}\end{center}
  
  \begin{figure}[!ht]
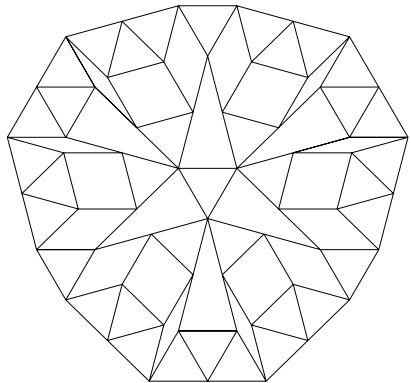

  	\centering
  	\caption{$(4;6)$-regular matchstick graph with 57 vertices and 117 edges.}
    %
  \end{figure}
  \noindent Discovered April 15, 2016 by M. Winkler. This rigid graph has a rotational symmetry of order 3 and contains three overlapped triplet kites.
  
  \newpage
  
  \begin{figure}[!ht]
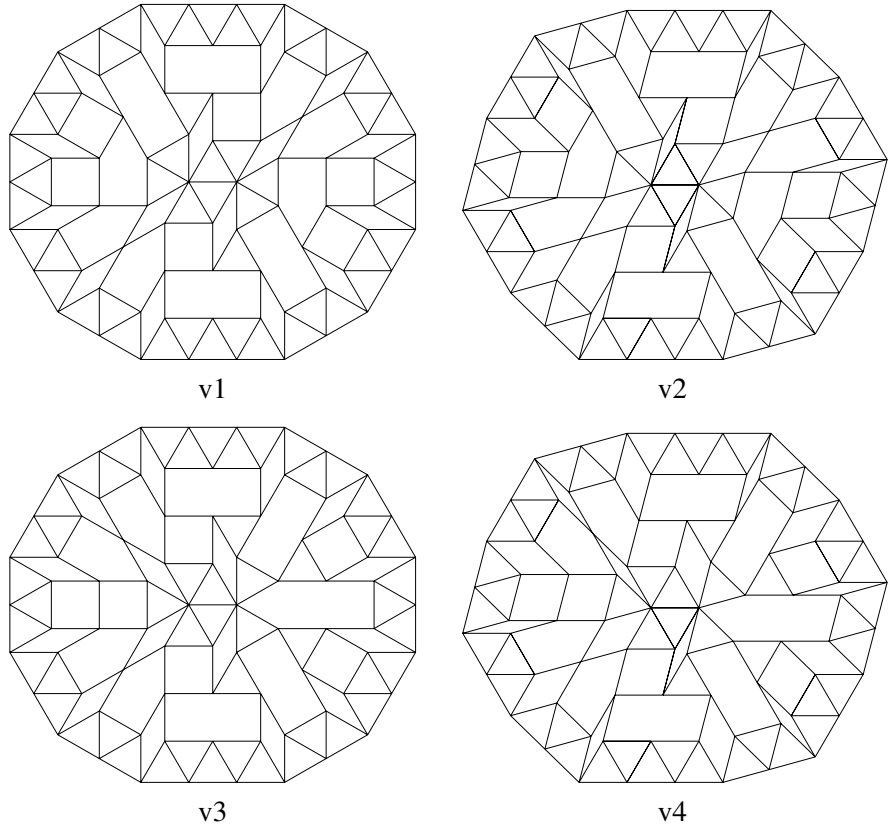

    \centering
    \caption{$(4;7)$-regular matchstick graphs with 78 vertices and 159 edges.}
    \begin{minipage}[t]{0.45\linewidth}
      \centering
          %
      \subcaption*{v1}
    \end{minipage}
    \quad
    \begin{minipage}[t]{0.45\linewidth}
      \centering
      %
      \subcaption*{v2}
    \end{minipage}
    
    %\begin{center}\end{center}
    
    \begin{minipage}[t]{0.45\linewidth}
      \centering
          %
      \subcaption*{v3}
    \end{minipage}
    \quad
    \begin{minipage}[t]{0.45\linewidth}
      \centering
      %
      \subcaption*{v4}
    \end{minipage}
  \end{figure}
  \noindent Discovered March 17, 2016 by M. Winkler, except for v3, which was discovered by P. Dinkelacker. The graph v3 is based on the graph v1 and has a slightly different internal geometry. These graphs are flexible. The graph v1 can be transformed into the graph v2, and the graph v3 can be transformed into the graph v4. The v1-based graphs always have a point symmetry. The graph v3 has a horizontal symmetry if the angles consist of 30 degrees and their multiples. For all other degrees the transformed graph is as asymmetric as the graph v4. The two vertices of degree 7 share an edge from which each graph gets its new minimality. The graphs v3 and v4 offer the possibility to rearrange two edges in the middle of the right side. It is unknown whether a rigid or kite-based $(4;7)$-regular matchstick graph with 159 edges or less exists.
  
  \newpage
  
  \begin{figure}[!ht]
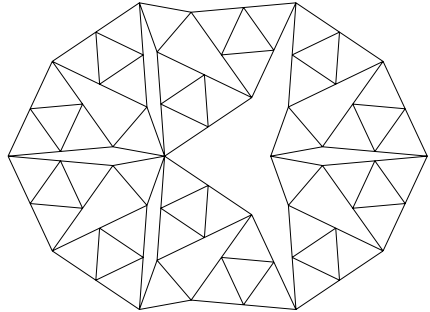

  	\centering
  	\caption{$(4;8)$-regular matchstick graph with 62 vertices and 126 edges.}

  \end{figure}
  \noindent Discovered March 14, 2016 by P. Dinkelacker. This rigid graph has a horizontal symmetry and consists only of kites. Two reverse double kites and one double kite.\\
  
  \begin{figure}[!ht]
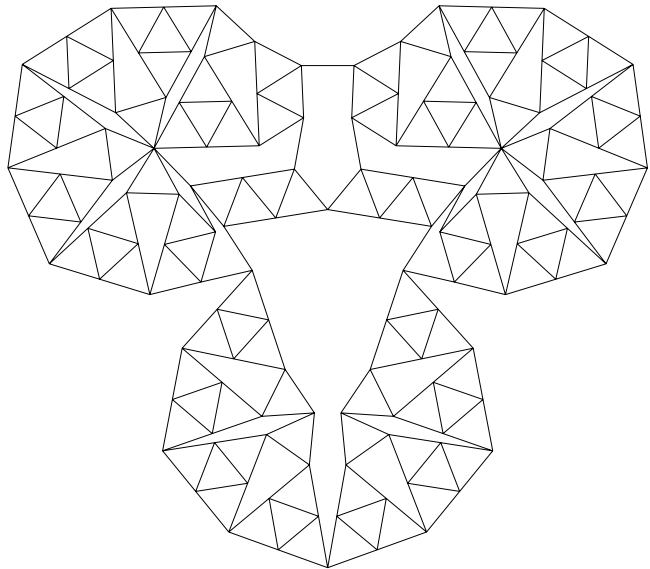

  	\centering
  	\caption{$(4;9)$-regular matchstick graph with 134 vertices and 273 edges.}
    %
  \end{figure}
  \noindent Discovered May 18, 2016 by P. Dinkelacker. This rigid graph has a vertical symmetry and contains four double kites and four slightly modified single kites.
  
  \newpage
  
  \begin{figure}[!ht]
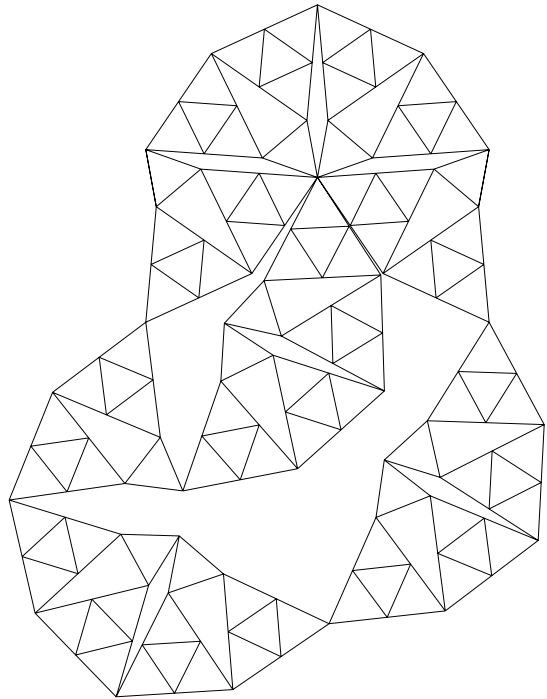

  	\centering
  	\caption{$(4;10)$-regular matchstick graph with 114 vertices and 231 edges.}
    %
  \end{figure}
  \noindent Discovered March 17, 2016 by P. Dinkelacker. This rigid graph is asymmetric and consists only of kites. Three double kites, two reverse double kites and one single kite. It remains an interesting question whether a symmetric $(4;10)$-regular matchstick graph with 231 edges or less exists.
  
  \newpage
  
  \begin{figure}[!ht]
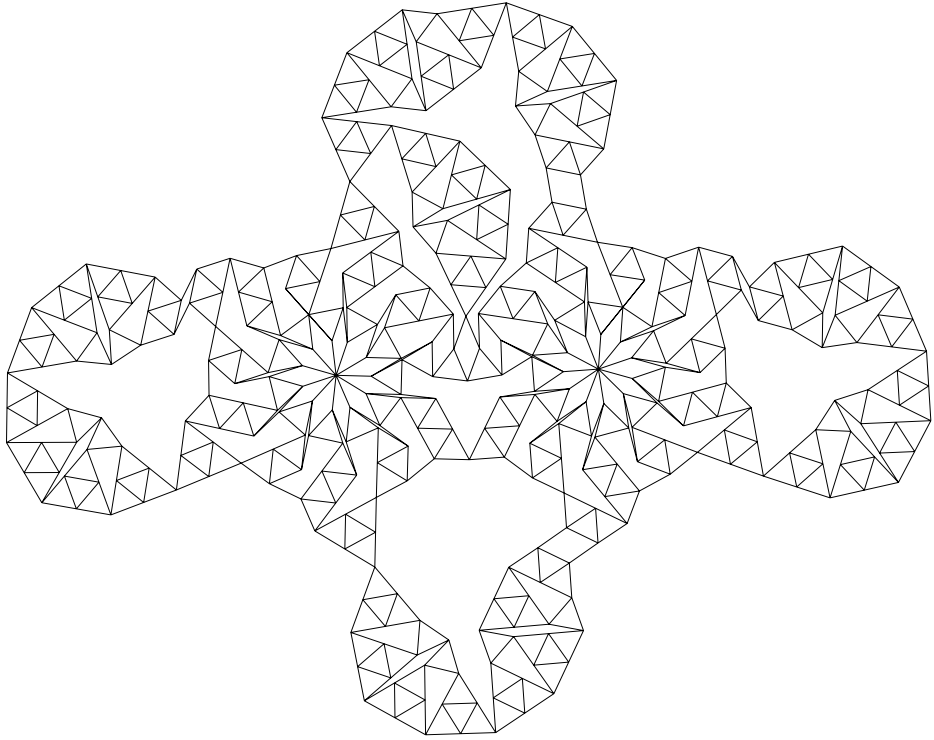

  	\centering
  	\caption{$(4;11)$-regular matchstick graph with 382 vertices and 771 edges.}
  	%
  \end{figure}
  \noindent Discovered April 15, 2018 by M. Winkler, S. Vogel and P. Dinkelacker. This rigid graph is asymmetric and contains four double kites, one reverse double kite, four kites and four slightly modified kites. There exists a few asymmetric variations of this graph with 771 edges, because the clasps can be varied. But the current design requires the least place in the plane.
  
  \newpage
  
  \noindent The interesting part of the graph in Figure 13 is the area around a vertex of degree 11, as the detail image in Figure 14 shows. This subgraph is rigid, but if removing only one edge the whole subgraph is flexible.
  
  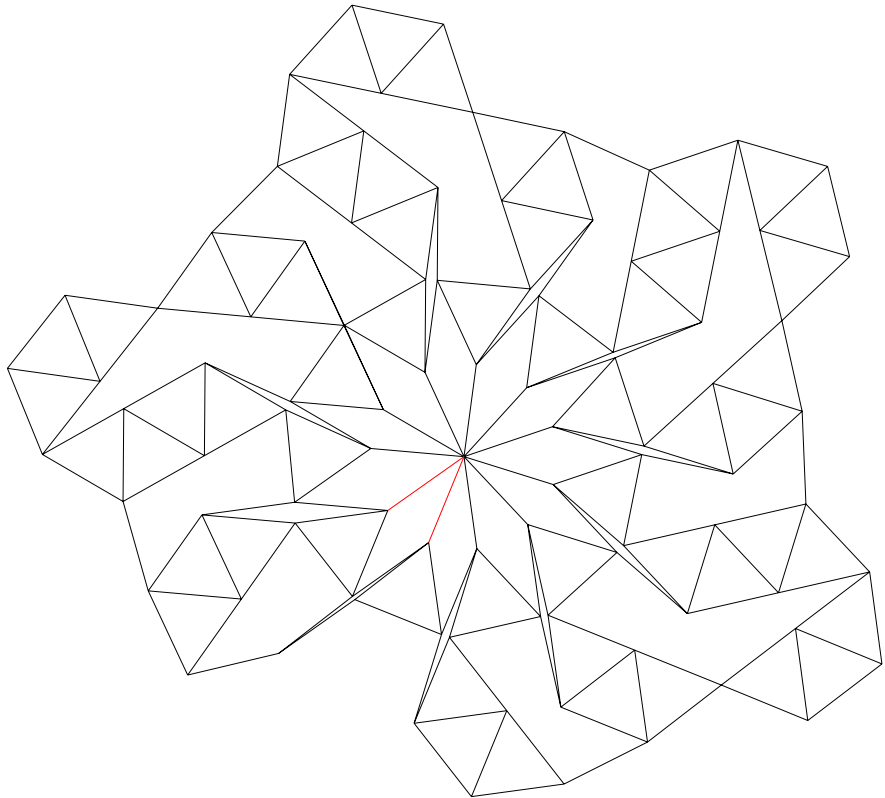
\begin{figure}[!ht]
  	\centering
  	\begin{tikzpicture}
    [y=0.8pt, x=0.8pt, yscale=-1.0, xscale=1.0]
    \draw[line width=0.01pt]
  	(168.2703,950.9935) -- (203.9330,925.6161)
  	(168.2703,950.9935) -- (184.8099,910.4684)
  	(141.4443,916.4072) -- (184.8099,910.4684)
  	(141.4443,916.4072) -- (168.2703,950.9935)
  	(176.8726,881.2331) -- (220.4726,885.0910)
  	(141.2099,906.6105) -- (184.8099,910.4684)
  	(141.2099,906.6105) -- (176.8726,881.2331)
  	(97.8443,912.5493) -- (141.4443,916.4072)
  	(97.8443,912.5493) -- (141.2099,906.6105)
  	(116.3033,952.2370) -- (141.4443,916.4072)
  	(116.3033,952.2370) -- (97.8443,912.5493)
  	(72.7033,948.3791) -- (97.8443,912.5493)
  	(72.7033,948.3791) -- (116.3033,952.2370)
  	(91.1623,988.0668) -- (116.3033,952.2370)
  	(91.1623,988.0668) -- (72.7033,948.3791)
  	(91.1623,988.0668) -- (133.7306,977.8795)
  	(137.0638,863.0370) -- (176.8726,881.2331)
  	(137.0638,863.0370) -- (141.2099,906.6105)
  	(60.9181,906.2252) -- (72.7033,948.3791)
  	(60.9181,906.2252) -- (98.9910,884.6311)
  	(60.9181,906.2252) -- (61.2535,862.4561)
  	(99.3263,840.8620) -- (137.0638,863.0370)
  	(98.9910,884.6311) -- (137.0638,863.0370)
  	(98.9910,884.6311) -- (99.3263,840.8620)
  	(61.2535,862.4561) -- (99.3263,840.8620)
  	(61.2535,862.4561) -- (98.9910,884.6311)
  	(23.1807,884.0502) -- (60.9181,906.2252)
  	(23.1807,884.0502) -- (61.2535,862.4561)
  	(23.1807,884.0502) -- (50.1474,849.5735)
  	(139.1352,859.0581) -- (176.8726,881.2331)
  	(139.1352,859.0581) -- (99.3263,840.8620)
  	(182.7352,862.9160) -- (220.4726,885.0910)
  	(182.7352,862.9160) -- (139.1352,859.0581)
  	(164.2762,823.2284) -- (139.1352,859.0581)
  	(164.2762,823.2284) -- (182.7352,862.9160)
  	(77.1140,815.0968) -- (120.6951,819.1626)
  	(77.1140,815.0968) -- (102.4256,779.3873)
  	(146.0067,783.4531) -- (164.2762,823.2284)
  	(120.6951,819.1626) -- (164.2762,823.2284)
  	(120.6951,819.1626) -- (146.0067,783.4531)
  	(102.4256,779.3873) -- (146.0067,783.4531)
  	(102.4256,779.3873) -- (120.6951,819.1626)
  	(50.1474,849.5735) -- (77.1140,815.0968)
  	(164.4657,823.1408) -- (182.7352,862.9160)
  	(164.4657,823.1408) -- (146.0067,783.4531)
  	(202.2031,845.3157) -- (220.4726,885.0910)
  	(202.2031,845.3157) -- (164.4657,823.1408)
  	(202.5385,801.5467) -- (164.4657,823.1408)
  	(202.5385,801.5467) -- (202.2031,845.3157)
  	(133.1332,748.1962) -- (102.4256,779.3873)
  	(133.1332,748.1962) -- (167.8359,774.8714)
  	(133.1332,748.1962) -- (173.5860,731.4804)
  	(208.2886,758.1557) -- (202.5385,801.5467)
  	(167.8359,774.8714) -- (202.5385,801.5467)
  	(167.8359,774.8714) -- (208.2886,758.1557)
  	(173.5860,731.4804) -- (208.2886,758.1557)
  	(173.5860,731.4804) -- (167.8359,774.8714)
  	(138.8833,704.8052) -- (133.1332,748.1962)
  	(138.8833,704.8052) -- (173.5860,731.4804)
  	(138.8833,704.8052) -- (181.7289,713.7549)
  	(207.9533,801.9247) -- (202.2031,845.3157)
  	(207.9533,801.9247) -- (208.2886,758.1557)
  	(226.2227,841.7000) -- (220.4726,885.0910)
  	(226.2227,841.7000) -- (207.9533,801.9247)
  	(251.5344,805.9905) -- (207.9533,801.9247)
  	(251.5344,805.9905) -- (226.2227,841.7000)
  	(224.5745,722.7046) -- (238.0545,764.3476)
  	(224.5745,722.7046) -- (267.3783,731.8521)
  	(280.8583,773.4951) -- (251.5344,805.9905)
  	(238.0545,764.3476) -- (251.5344,805.9905)
  	(238.0545,764.3476) -- (280.8583,773.4951)
  	(267.3783,731.8521) -- (280.8583,773.4951)
  	(267.3783,731.8521) -- (238.0545,764.3476)
  	(181.7289,713.7549) -- (224.5745,722.7046)
  	(255.5466,809.2045) -- (226.2227,841.7000)
  	(255.5466,809.2045) -- (280.8583,773.4951)
  	(249.7965,852.5956) -- (220.4726,885.0910)
  	(249.7965,852.5956) -- (255.5466,809.2045)
  	(290.2493,835.8798) -- (255.5466,809.2045)
  	(290.2493,835.8798) -- (249.7965,852.5956)
  	(307.2102,749.9979) -- (267.3783,731.8521)
  	(307.2102,749.9979) -- (298.7297,792.9388)
  	(307.2102,749.9979) -- (340.1579,778.8126)
  	(331.6775,821.7536) -- (290.2493,835.8798)
  	(298.7297,792.9388) -- (290.2493,835.8798)
  	(298.7297,792.9388) -- (331.6775,821.7536)
  	(340.1579,778.8126) -- (331.6775,821.7536)
  	(340.1579,778.8126) -- (298.7297,792.9388)
  	(348.6384,735.8717) -- (307.2102,749.9979)
  	(348.6384,735.8717) -- (340.1579,778.8126)
  	(348.6384,735.8717) -- (358.9354,778.4136)
  	(348.6384,735.8717) -- (390.6293,748.2251)
  	(291.2247,838.4694) -- (249.7965,852.5956)
  	(291.2247,838.4694) -- (331.6775,821.7536)
  	(261.9008,870.9648) -- (220.4726,885.0910)
  	(261.9008,870.9648) -- (291.2247,838.4694)
  	(304.7046,880.1123) -- (291.2247,838.4694)
  	(304.7046,880.1123) -- (261.9008,870.9648)
  	(369.2325,820.9555) -- (336.9686,850.5339)
  	(369.2325,820.9555) -- (378.7162,863.6861)
  	(346.4522,893.2645) -- (304.7046,880.1123)
  	(336.9686,850.5339) -- (304.7046,880.1123)
  	(336.9686,850.5339) -- (346.4522,893.2645)
  	(378.7162,863.6861) -- (346.4522,893.2645)
  	(378.7162,863.6861) -- (336.9686,850.5339)
  	(400.9263,790.7670) -- (369.2325,820.9555)
  	(358.9354,778.4136) -- (369.2325,820.9555)
  	(358.9354,778.4136) -- (400.9263,790.7670)
  	(390.6293,748.2251) -- (400.9263,790.7670)
  	(390.6293,748.2251) -- (358.9354,778.4136)
  	(303.6484,884.1170) -- (261.9008,870.9648)
  	(303.6484,884.1170) -- (346.4522,893.2645)
  	(262.2202,898.2432) -- (220.4726,885.0910)
  	(262.2202,898.2432) -- (303.6484,884.1170)
  	(295.1680,927.0579) -- (303.6484,884.1170)
  	(295.1680,927.0579) -- (262.2202,898.2432)
  	(380.4778,907.4210) -- (378.7162,863.6861)
  	(380.4778,907.4210) -- (337.8229,917.2395)
  	(380.4778,907.4210) -- (367.6534,949.2705)
  	(324.9985,959.0889) -- (295.1680,927.0579)
  	(337.8229,917.2395) -- (295.1680,927.0579)
  	(337.8229,917.2395) -- (324.9985,959.0889)
  	(367.6534,949.2705) -- (324.9985,959.0889)
  	(367.6534,949.2705) -- (337.8229,917.2395)
  	(410.3083,939.4520) -- (380.4778,907.4210)
  	(410.3083,939.4520) -- (367.6534,949.2705)
  	(410.3083,939.4520) -- (375.6013,966.1216)
  	(410.3083,939.4520) -- (416.0513,982.8439)
  	(292.0507,930.2742) -- (262.2202,898.2432)
  	(292.0507,930.2742) -- (324.9985,959.0889)
  	(250.3031,917.1220) -- (220.4726,885.0910)
  	(250.3031,917.1220) -- (292.0507,930.2742)
  	(259.7868,959.8526) -- (292.0507,930.2742)
  	(259.7868,959.8526) -- (250.3031,917.1220)
  	(340.8943,992.7912) -- (300.3405,976.3219)
  	(340.8943,992.7912) -- (306.3546,1019.6771)
  	(265.8008,1003.2078) -- (259.7868,959.8526)
  	(300.3405,976.3219) -- (259.7868,959.8526)
  	(300.3405,976.3219) -- (265.8008,1003.2078)
  	(306.3546,1019.6771) -- (265.8008,1003.2078)
  	(306.3546,1019.6771) -- (300.3405,976.3219)
  	(306.3546,1019.6771) -- (267.3142,1039.4683)
  	(381.3443,1009.5135) -- (340.8943,992.7912)
  	(375.6013,966.1216) -- (340.8943,992.7912)
  	(375.6013,966.1216) -- (381.3443,1009.5135)
  	(416.0513,982.8439) -- (381.3443,1009.5135)
  	(416.0513,982.8439) -- (375.6013,966.1216)
  	(256.3172,960.4772) -- (250.3031,917.1220)
  	(256.3172,960.4772) -- (265.8008,1003.2078)
  	(226.4867,928.4462) -- (220.4726,885.0910)
  	(226.4867,928.4462) -- (256.3172,960.4772)
  	(213.6623,970.2957) -- (256.3172,960.4772)
  	(213.6623,970.2957) -- (226.4867,928.4462)
  	(209.9470,968.9713) -- (203.9330,925.6161)
  	(209.9470,968.9713) -- (226.4867,928.4462)
  	(197.1226,1010.8208) -- (213.6623,970.2957)
  	(197.1226,1010.8208) -- (209.9470,968.9713)
  	(267.3142,1039.4683) -- (240.4882,1004.8820)
  	(267.3142,1039.4683) -- (223.9486,1045.4071)
  	(240.4882,1004.8820) -- (213.6623,970.2957)
  	(240.4882,1004.8820) -- (197.1226,1010.8208)
  	(223.9486,1045.4071) -- (197.1226,1010.8208)
  	(223.9486,1045.4071) -- (240.4882,1004.8820)
  	(169.3933,952.5020) -- (203.9330,925.6161)
  	(169.3933,952.5020) -- (209.9470,968.9713)
  	(133.7306,977.8795) -- (168.2703,950.9935)
  	(133.7306,977.8795) -- (169.3933,952.5020)
  	(33.7730,808.9813) -- (6.8063,843.4580)
  	(6.8063,843.4580) -- (50.1474,849.5735)
  	(50.1474,849.5735) -- (33.7730,808.9813)
  	(210.9024,681.1244) -- (168.0568,672.1746)
  	(138.8833,704.8052) -- (168.0568,672.1746)
  	(168.0568,672.1746) -- (181.7289,713.7549)
  	(181.7289,713.7549) -- (210.9024,681.1244)
  	(210.9024,681.1244) -- (224.5745,722.7046)
  	(23.1807,884.0502) -- (6.8063,843.4580)
  	(33.7730,808.9813) -- (77.1140,815.0968);
  	% red lines
  	\draw[red, line width=0.01pt]
  	(203.9330,925.6161) -- (220.4726,885.0910)
  	(184.8099,910.4684) -- (220.4726,885.0910);
  	
  	\end{tikzpicture}
  	\caption{Detail around the right vertex of degree 11 in Figure 13.}
  \end{figure}
  
  \noindent This flexibility makes it possible to adjust the eleven angles around the centered vertex so that all edges have exactly unit length. Beginning clockwise with the angle between the red edges, these degrees are\\
  \\
  32.362519660072210,\quad 40.49207000332465,\quad 25.382433534610843,\\
  34.890820876760450,\quad 32.21894760945070,\quad 34.514335947363630,\\
  29.108515978283318,\quad 36.31491131809427,\quad 29.550687898877964,\\
  35.065359484316880,\quad 30.09939768884507.
  
  \newpage
  
  \section{\large{Further new minimal $(4;n)$-regular matchstick graphs for $n=5$ and $n=6$}}
  
  Among the smallest known $(4;n)$-regular matchstick graphs of Chapter 4, the authors had discovered further new minimal graphs for $n=5$ and $n=6$ consisting of 118, 120, 121 and 126 edges. Most of the graphs in Chapter 5 offer the possibility to construct examples with 121, 122, 123, 124 and 125 edges by inserting additional edges.
  \\
  \begin{figure}[!ht]
  	\centering
  	\caption{$(4;5)$-regular matchstick graphs with 58 vertices and 118 edges, and 59 vertices and 120 edges.}
  	\begin{minipage}[t]{0.45\linewidth}
  		\centering

  		\subcaption{118 edges}
  	\end{minipage}
  	\quad
  	\begin{minipage}[t]{0.45\linewidth}
  		\centering
  		%
  		\subcaption{120 edges}
  	\end{minipage}
  \end{figure}
  
  \noindent Discovered February 23, 2017 by M. Winkler. The graph with 118 edges is rigid, has a vertical symmetry and contains a kite and a triplet kite. The graph with 120 edges is rigid, asymmetric and contains a kite and a slightly modified triplet kite. The graph 15b can be constructed from the graph 15a by adding one vertex and two edges and replacing two edges.
  
  \newpage
  
  \begin{figure}[!ht]
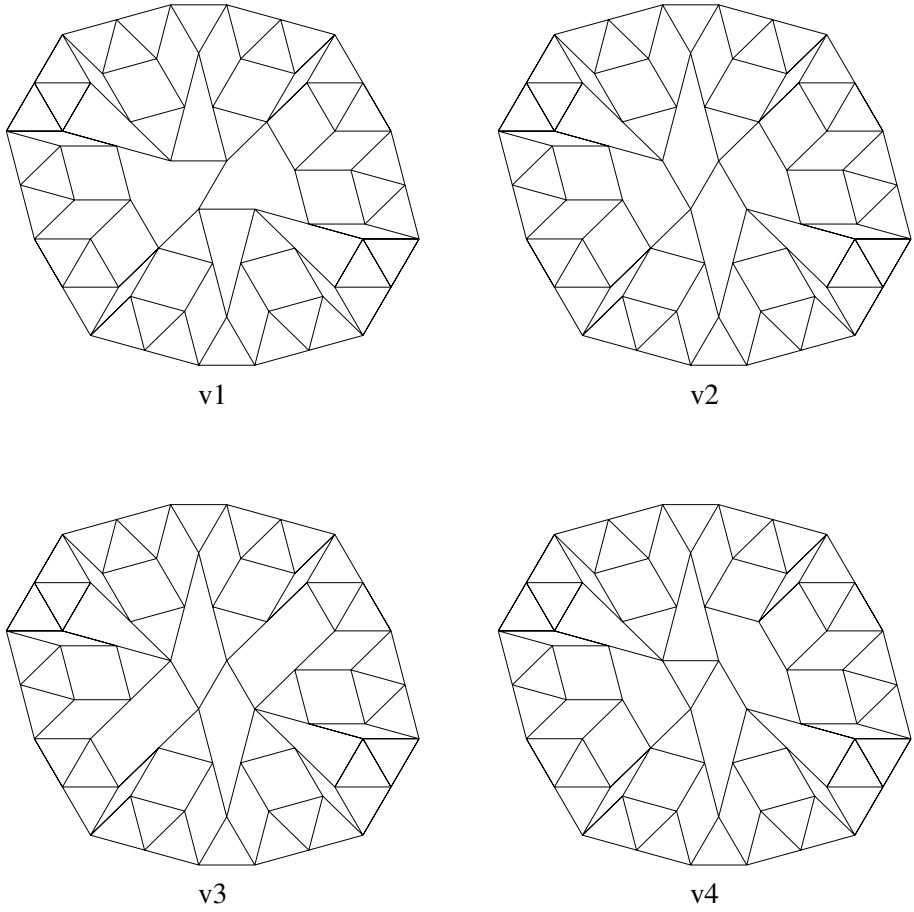

    \centering
    \caption{$(4;5)$-regular matchstick graphs with 60 vertices and 121 edges.}
    \begin{minipage}[t]{0.45\linewidth}
      \centering
      %
      \subcaption*{v1}
    \end{minipage}
    \quad\quad
    \begin{minipage}[t]{0.45\linewidth}
      \centering
      %
      \subcaption*{v2}
    \end{minipage}
    
    \begin{center}\end{center}
    
    \begin{minipage}[t]{0.45\linewidth}
       \centering
       %
      \subcaption*{v3}
    \end{minipage}
    \quad\quad
    \begin{minipage}[t]{0.45\linewidth}
      \centering
      %
      \subcaption*{v4}
    \end{minipage}
  \end{figure}
  \noindent Discovered April 13, 2016 by  M. Winkler. The rigid graphs v1, v2 and v3 have a point symmetry. The rigid graph v4 is asymmetric.
  
  \newpage
  
  \begin{figure}[!ht]
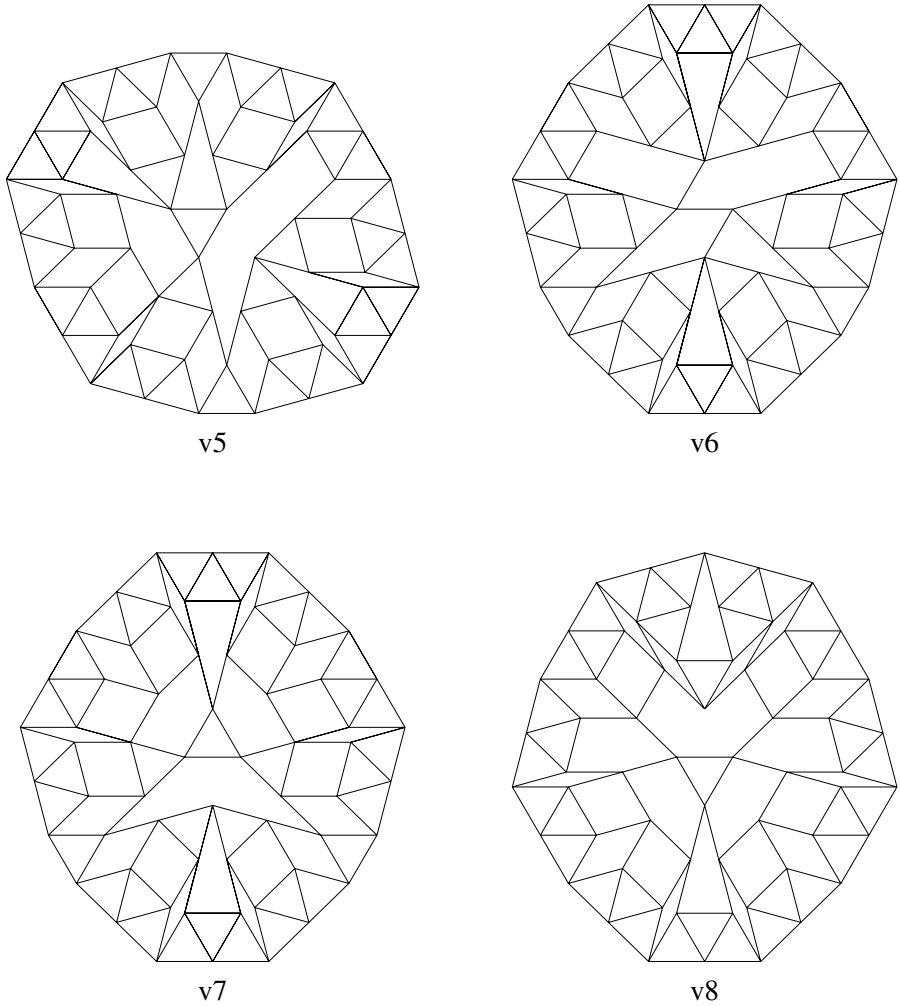

    \centering
    \caption{$(4;5)$-regular matchstick graphs with 60 vertices and 121 edges.}
    \begin{minipage}[t]{0.45\linewidth}
      \centering
      %
      \subcaption*{v5}
    \end{minipage}
    \quad\quad
    \begin{minipage}[t]{0.45\linewidth}
    	\centering
    	%
    	\subcaption*{v6}
    \end{minipage}
    
    \begin{center}\end{center}
    
    \begin{minipage}[t]{0.45\linewidth}
       \centering
       %
      \subcaption*{v8}
    \end{minipage}
  \end{figure}
  \noindent Discovered April 13-14, 2016 and May 4, 2017 (v8) by M. Winkler. The rigid graphs v5 and v6 are asymmetric. The rigid graphs v7 and v8 have a vertical symmetry.
  
  \newpage
  
  \begin{figure}[!ht]
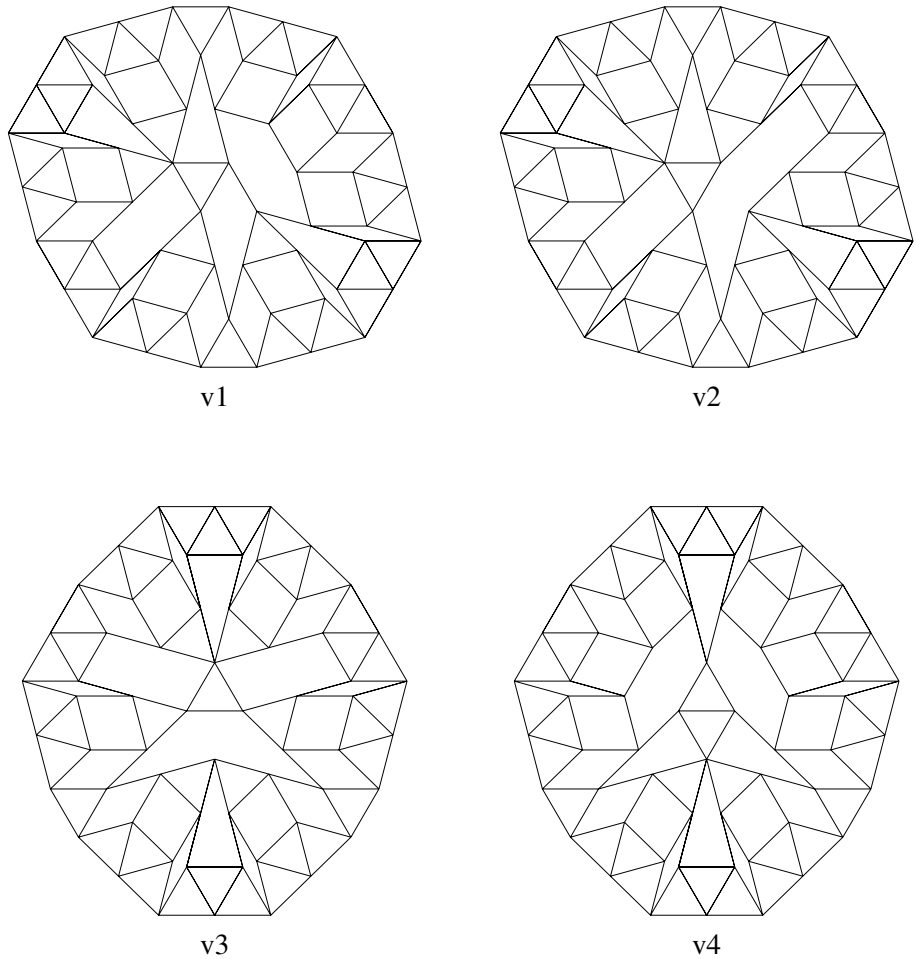

    \centering
    \caption{$(4;6)$-regular matchstick graphs with 60 vertices and 121 edges.}
    \begin{minipage}[t]{0.45\linewidth}
      \centering
      %
      \subcaption*{v2}
    \end{minipage}
    
    \begin{center}\end{center}
    
    \begin{minipage}[t]{0.45\linewidth}
       \centering
       %
      \subcaption*{v4}
    \end{minipage}
  \end{figure}
  \noindent Discovered April 14, 2016 by M. Winkler, except for v3, which was discovered by P. Dinkelacker. The rigid graphs v1 and v2 are asymmetric. The rigid graphs v3 and v4 have a vertical symmetry. The vertical symmetry of these graphs was discovered by P. Dinkelacker and is based on the equivalent asymmetric and point symmetric graphs of this Chapter by M. Winkler.
  
  \newpage
  
  \begin{figure}[!ht]
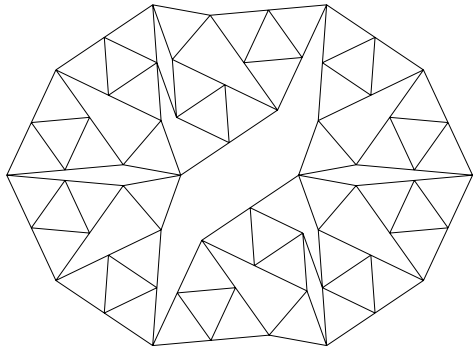

  	\centering
  	\caption{$(4;5)$-regular matchstick graph with 62 vertices and 126 edges.}
    %
  \end{figure}
  \noindent Discovered March 14, 2016 by M. Winkler. This rigid graph has a point symmetry and consists of two slightly modified reverse double kites and two single kites.
  
  \begin{center}\end{center}
  
  \begin{figure}[!ht]
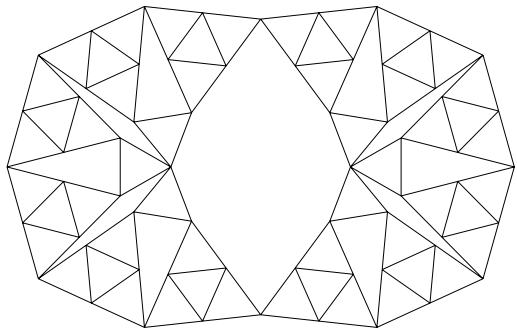

  	\centering
  	\caption{$(4;6)$-regular matchstick graph with 62 vertices and 126 edges.}
    %
  \end{figure}
  \noindent Discovered March 14, 2016 by P. Dinkelacker. This rigid graph has a vertical and a horizontal symmetry and consists only of kites. Two double kites and two single kites.
  
  \newpage
  
  \begin{figure}[!ht]
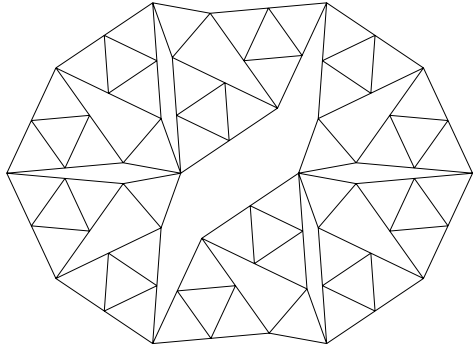

  	\centering
  	\caption{$(4;6)$-regular matchstick graph with 62 vertices and 126 edges.}
    %
  \end{figure}
  \noindent Discovered March 14, 2016 by M. Winkler. This rigid graph has a point symmetry and consists only of kites. Two reverse double kites and two single kites. This graph needs less space in the plane than the $(4;6)$-graph in Figure 20. The only difference between the graphs in Figure 19 and 21 are two replaced edges.
  
  \section{\large{Infinite $(4;n)$-regular matchstick graphs for $n>11$}}
  
  \noindent For each $n>11$ only infinite $(4;n)$-regular matchstick graphs with an infinite number of vertices exist. The graphs in Figure 22 show only the simplest version of these kind of graphs with one vertex of degree $n$ in the center. But there exists an infinite number of different constructions with up to an infinite number of vertices of degree $n$. All these kinds of infinite graphs are high flexible. The next four examples of infinite $(4;n)$-regular matchstick graphs for $n=12$ and $n=13$ will illustrate this.
  
  \newpage
  
  
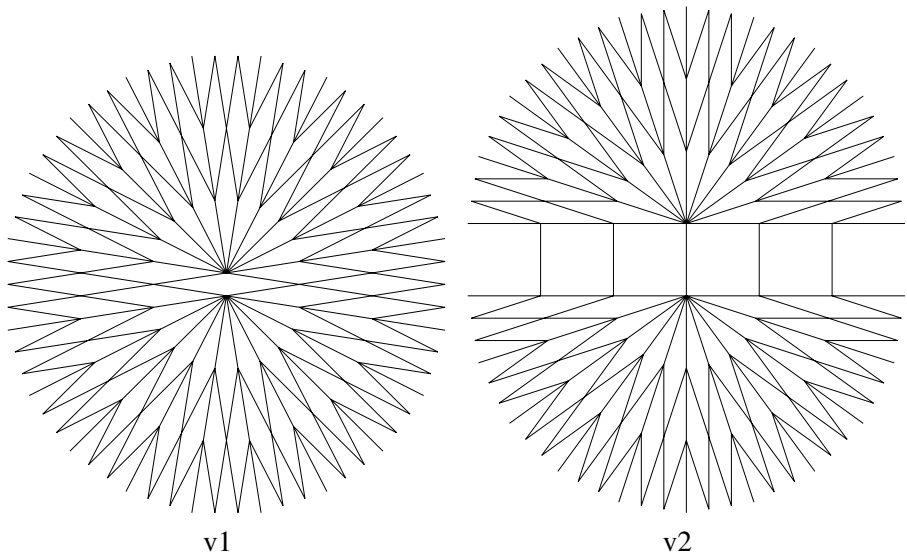
\begin{figure}[!ht]
    \centering
    \caption{Infinite $(4;12)$-regular and $(4;13)$-regular matchstick graphs.}
    \begin{minipage}[t]{0.45\linewidth}
      \centering

      \subcaption{$(4;12)$}
    \end{minipage}
    \quad
    \begin{minipage}[t]{0.45\linewidth}
      \centering
      %
      \subcaption{$(4;13)$}
    \end{minipage}
  \end{figure}
  
  \begin{center}\end{center}
  
  \begin{figure}[!ht]
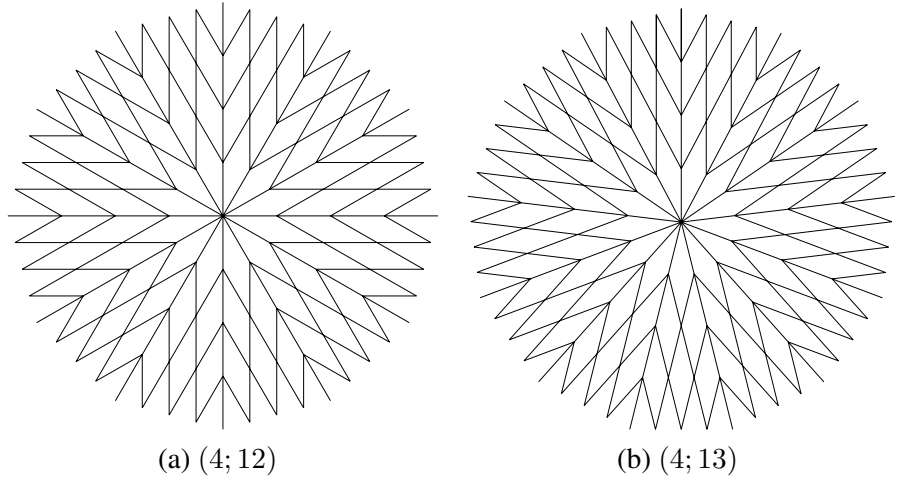

  	\centering
  	\caption{Infinite $(4;12)$-regular matchstick graphs with 2 vertices of degree 12.}
  	\begin{minipage}[t]{0.45\linewidth}
  		\centering
  		%
 
  		\subcaption*{v1}
  	\end{minipage}
  	\quad
  	\begin{minipage}[t]{0.45\linewidth}
  		\centering
  		%
  		\subcaption*{v2}
  	\end{minipage}
  \end{figure}
  
  \newpage
  
  \begin{figure}[!ht]
  	\centering
  	\caption{Infinite $(4;12)$-regular matchstick graph flexibility example.}
  	%
    	
  \end{figure}
  
  \newpage
  
  \begin{figure}[!ht]
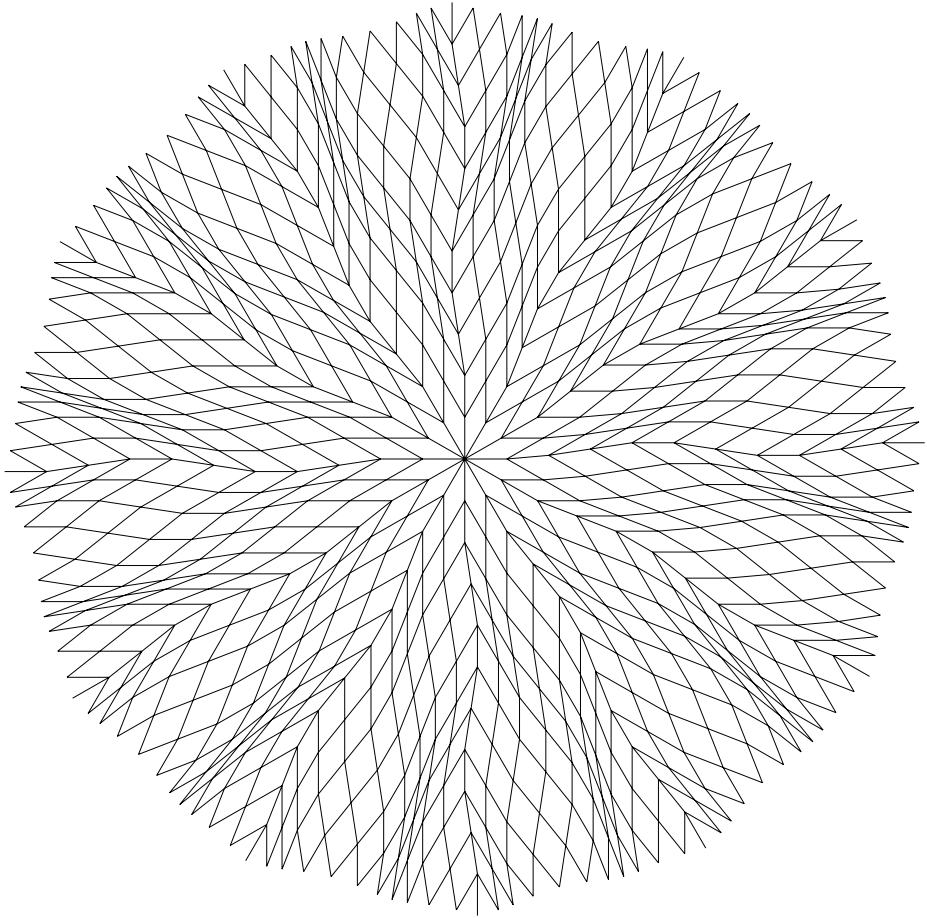

  	\centering
  	\caption{Infinite $(4;12)$-regular matchstick graph flexibility example.}
  	%
    	
  \end{figure}
  \noindent The graph in Figure 25 is a transformed variation of the graph in Figure 24.
  
  \newpage
  
  \section{\large{References}}
  
  1. Erich Friedman, \textit{Math Magic, Problem of the Month (December 2005)}, Problem 4, Smallest Known m/n Matchstick Graphs.\\
  \footnotesize(http://www2.stetson.edu/~efriedma/mathmagic/1205.html)\normalsize
  \\ \\
  2. Heiko Harborth, \textit{Match Sticks in the Plane}, The Lighter Side of Mathematics. Proceedings of the Eug\`ene Strens Memorial Conference of Recreational Mathematics \& its History, Calgary, Canada, July 27 -- August 2, 1986 (Washington) (Richard K. Guy and Robert E. Woodrow, eds.), Spectrum Series, The Mathematical Association of America, 1994, pp. 281--288.
  \\ \\
  3. Sascha Kurz and Giuseppe Mazzuoccolo, \textit{3-regular matchstick graphs with given girth}, Geombinatorics Quarterly Volume 19, Issue 4, April 2010, pp. 156--175.\\
  \footnotesize(http://arxiv.org/pdf/1401.4360v1.pdf)\normalsize
  \\ \\
  4. Wayback Machine Internet Archive.\\
  \footnotesize(https://web.archive.org/web/20151209031635/http://www2.stetson.edu/$\sim$efriedma/\\mathmagic/1205.html)\normalsize
  \\ \\
  5. Stefan Vogel, \textit{Beweglichkeit eines Streichholzgraphen bestimmen}, July 2016.\\
  \footnotesize(https://tinyurl.com/yc8at6r7)\normalsize
  \\ \\
  6. Stefan Vogel, \textit{Matchstick Graphs Calculator (MGC)}, a software for the construction and calculation of matchstick graphs, 2016 -- 2017.\\
  \footnotesize(http://mikewinkler.co.nf/matchstick\_graphs\_calculator.htm)\normalsize
  \\ \\
  7. Mike Winkler, Peter Dinkelacker, and Stefan Vogel, \textit{Streichholzgraphen 4-regul\"ar und 4/n-regul\"ar (n$>$4) und 2/5}, thread in a graph theory internet forum, used nicknames: P. Dinkelacker (haribo), M. Winkler (Slash),\\
  \footnotesize(https://tinyurl.com/ya3g6p7w)\normalsize
  \\ \\
  8. Mike Winkler, \textit{Der gro\ss{}e Bruder des Harborth-Graphen}, July 2016.\\
  \footnotesize(https://tinyurl.com/ydbawq34)\normalsize
  \\ \\
  9. Mike Winkler, \textit{Ein 4-regul\"arer Streichholzgraph mit 114 Kanten}, May 2016.\\
  \footnotesize(https://tinyurl.com/ydcwqrlx)\normalsize
  \\ \\
  10. Mike Winkler, \textit{Ein neuer 4-regul\"arer Streichholzgraph}, Mitteilungen der Deutschen Mathematiker-Vereinigung (DMV), Band 24, Heft 2, Seiten 74--75, DOI: 10.1515/dmvm-2016-0031, July 2016.\\
  \footnotesize(http://www.degruyter.com/view/j/dmvm.2016.24.issue-2/dmvm-2016-0031/dmvm-2016-0031.xml)\normalsize
  
  \begin{center}\end{center}
  
  \begin{center}
    \textbf{Acknowledgements}
  \end{center}
  \noindent
  The authors wishes to express their thanks to Erich Friedman for this interesting and challenging task and for presenting the new minimal graphs on his \textit{Math Magic} website [1].

\end{document}